\documentclass[twocolumn]{autart}

\usepackage[noframe]{showframe}
\usepackage{graphicx}
\usepackage{amsmath, amssymb}
\usepackage{algpseudocode}
\algnewcommand{\LineComment}[1]{\State {\color{violet}\# #1}}
\usepackage{algorithm}

\usepackage{dsfont}
\newcommand{\indicator}{\mathds{1}}

\usepackage{caption}
\usepackage{subcaption}

\usepackage{xcolor}

\usepackage{tikz}
\usetikzlibrary{shapes,arrows,positioning,calc}

\usepackage[colorlinks=true]{hyperref}

\newcommand{\shopt}{\textsc{SH2OPT}}
\newcommand{\muSpace}{\mathcal{D}}
\newcommand{\muSpaceStab}{\mathcal{D}_{\mathrm{stab}}}
\newcommand{\stablereal}{\mathcal{RH}_{\infty}}
\newcommand{\norm}[1]{\left\| #1 \right\|}
\newcommand{\hinner}[2]{{\left\langle #1, #2 \right\rangle}}
\newcommand{\hnorm}[1]{\norm{#1}_{\mathcal{H}_2}}
\newcommand{\dGdmu}[1]{\frac{\partial G}{\partial \mu_{#1}}}
\DeclareMathOperator{\tr}{tr}
\DeclareMathOperator{\real}{Re}

\newcommand{\tp}[1]{{#1}^{\mathrm{T}}}

\edef\endfrontmatter{%
  \unexpanded\expandafter{\endfrontmatter}
  \noexpand\endNoHyper 
}

\begin{document}

\tikzset{
block/.style = {draw, fill=white, rectangle, minimum height=3em, minimum width=3em},
tmp/.style  = {coordinate}, 
sum/.style= {draw, fill=white, circle, node distance=1cm, inner sep=0pt},
input/.style = {coordinate},
output/.style= {coordinate},
pinstyle/.style = {pin edge={to-,thin,black}}
}

\begin{frontmatter}

\title{Stochastic Optimization of Large-Scale Parametrized Dynamical Systems}

\author[Eindhoven]{Pascal Den Boef}\ead{p.d.boef@tue.nl},
\author[Eindhoven]{Jos Maubach},
\author[Eindhoven]{Wil Schilders},
\author[Eindhoven]{Nathan van de Wouw},
\address[Eindhoven]{Eindhoven University of Technology, Eindhoven}                                          
      

\begin{abstract}
Many relevant problems in the area of systems and control, such as controller synthesis, observer design and model reduction, can be viewed as optimization problems involving dynamical systems: for instance, maximizing performance in the synthesis setting or minimizing error in the reduction setting. When the involved dynamics are large-scale (e.g., high-dimensional semi-discretizations of partial differential equations), the optimization becomes computationally infeasible. Existing methods in literature lack computational scalability or solve an approximation of the problem (thereby losing guarantees with respect to the original problem). In this paper, we propose a novel method that circumvents these issues. The method is an extension of Stochastic Gradient Descent (SGD) which is widely used in the context of large-scale machine learning problems. The proposed SGD scheme minimizes the $\mathcal{H}_2$ norm of a (differentiable) parametrized dynamical system, and we prove that the scheme is guaranteed to preserve stability with high probability under boundedness conditions on the step size. Conditioned on the stability preservation, we also obtain probabilistic convergence guarantees to local minimizers. The method is also applicable to problems involving non-realizable dynamics as it only requires frequency-domain input-output samples. We demonstrate the potential of the approach on two numerical examples: fixed-order observer design for a large-scale thermal model and controller tuning for an infinite-dimensional system.
\end{abstract}

\end{frontmatter}

\section{Introduction}

\begin{figure}[h!]
	\centering
	\begin{subfigure}[b]{0.39\textwidth}
		\centering
		\includegraphics[width=\textwidth]{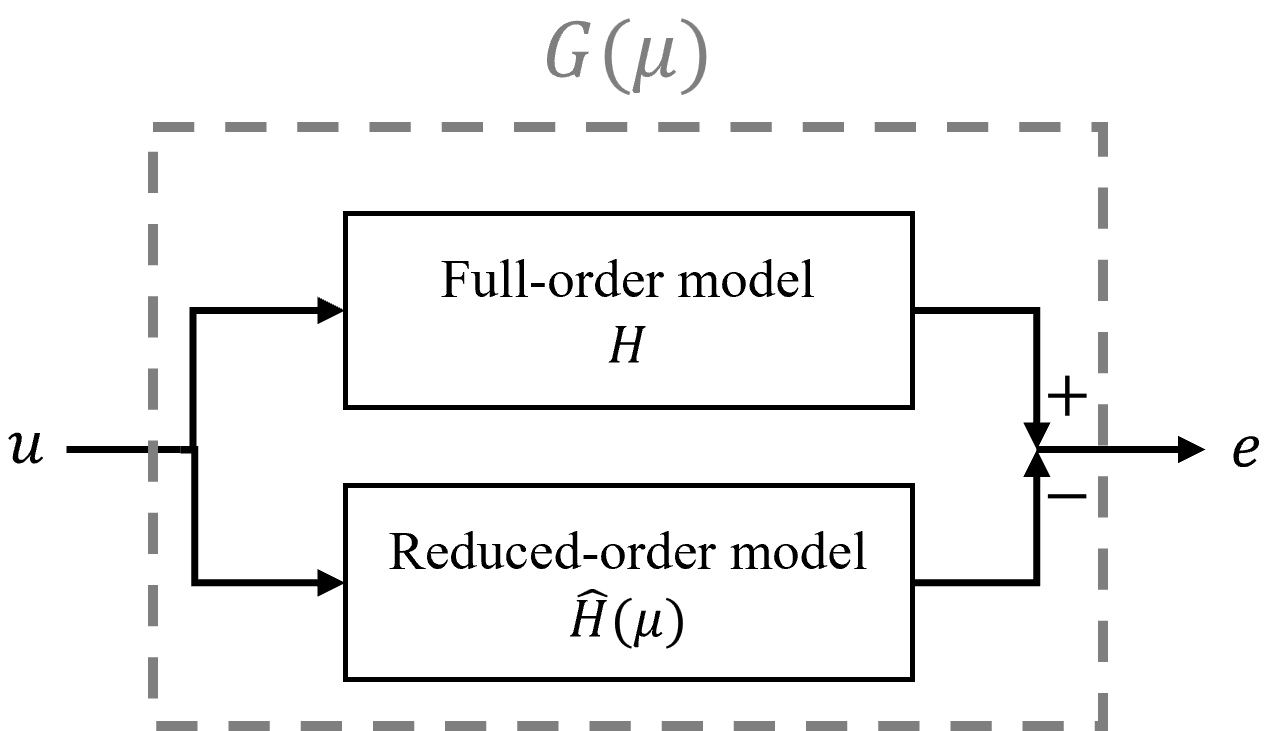}
		\caption{\textbf{Model reduction}: $G(\mu)$ represents the reduction error dynamics of a reduced-order model $\widehat{H}(\mu)$ approximating the input-output behaviour of $H$.}
		\label{fig:optimization_examples:reduction}	
	\end{subfigure}
	\hfill
	\begin{subfigure}[b]{0.39\textwidth}
		\centering
		\includegraphics[width=\textwidth]{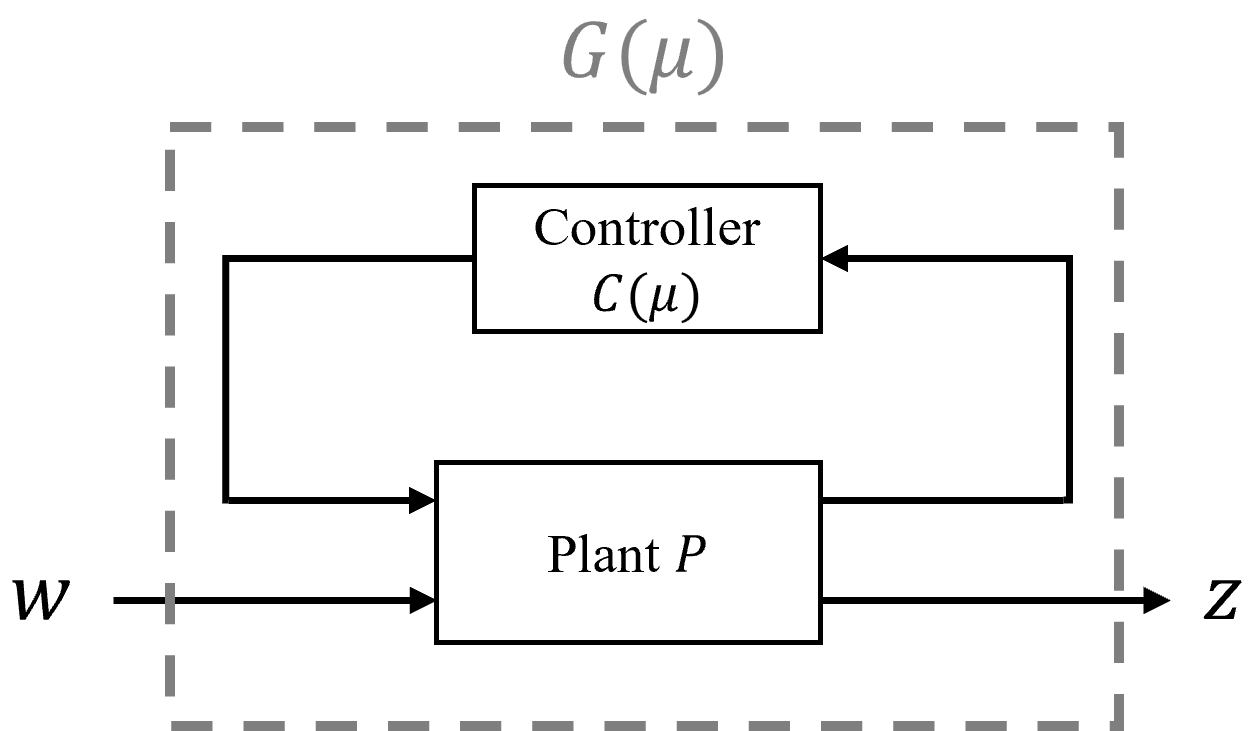}
		\caption{\textbf{Fixed-order control design}: $G(\mu)$ represents the control performance of a fixed-order controller for plant $P$.}
		\label{fig:optimization_examples:synthesis}
	\end{subfigure}
	\caption{A unified view of several system-theoretic problems involving an objective related to a parametrized input-output model $G(\mu)$. The vector $\mu$ parametrises degrees-of-freedom of the dynamical blocks, such as state-space coefficients or pole/zero locations.}
	\label{fig:optimization_examples}
\end{figure}

Problems involving the optimization of dynamical systems are encountered in system-theoretic problems related to model reduction \cite{Hokanson2020,Beattie2009,Necoara2022}, controller design \cite{Petersson2013,Lewis2012,peretz2016randomized,zeng2017structured,hyland1984optimal} and state estimation \cite{bernstein1985optimal}. In this scope, the optimization is typically posed in terms of input-to-output behaviour characterized by a transfer function $G(\mu)$. Herein, $G(\mu)$ represents a parametrized dynamical system, with $\mu$ a parameter vector denoting application-dependent degrees-of-freedom. Figure \ref{fig:optimization_examples} demonstrates how $G(\mu)$ provides a unified view of several of the mentioned system-theoretic applications. In this unified view, the system $G(\mu)$ can be seen as an interconnected system consisting of a mixture of parametric and non-parametric components. For example, the fixed-order control synthesis problem displayed in Figure \ref{fig:optimization_examples:synthesis} consists of a given generalized plant model $P$ (non-parametric) and the fixed-order controller to-be-designed $C(\mu)$ (parametric). 
When the average performance over frequencies is of interest, the $\mathcal{H}_2$ norm of $G(\mu)$ is a natural candidate for the optimality criterion leading to the following cost function $c(\mu)$:
\begin{equation}
	c(\mu) := \frac{1}{2} \hnorm{G(\mu)}^2.
	\label{eqn:cost}
\end{equation}
Finding solutions $\mu$ that (locally) minimize (\ref{eqn:cost}) has received considerable attention in literature. Most methods can be derived from the necessary optimality condition that the gradient vanishes at a local minimizer. For example, the method of optimal projection equations in control design \cite{hyland1984optimal} and state estimation \cite{bernstein1985optimal} characterizes local minimizers by a set of coupled Riccati equations. To solve these equations one can use continuation methods \cite{collins1998comparison}. Other gradient-based approaches for solving structured $\mathcal{H}_2$ optimization problems are proposed in \cite{Petersson2013,Necoara2020,Beattie2009}. The aforementioned methods rely on solutions of matrix equations involving the state-space matrices of the system-to-be-optimized. If these are large (as is the case with, e.g., high-dimensional Finite Element Method (FEM) models), solving these equations becomes computationally infeasible. The development of efficient (but approximate) solvers for large-scale sparse Lyapunov equations \cite{Kurschner2016,haber2016sparse} and Riccati equations \cite{benner2004solving,benner2020numerical}, enables working with approximations of solutions to the matrix equations. Whereas the aforementioned methods attempt to solve the \emph{exact} problem and use \emph{approximate} solvers for scaling to high-dimensional dynamics, another class of methods can be classified that \emph{approximate} the problem such that an \emph{exact} solver can be used. For example, when designing a controller for a high-dimensional plant, one possibility is to first reduce the order of the plant using a model reduction technique and to subsequently design an $\mathcal{H}_2$ controller for the reduced-order plant \cite{anderson1993controller}. Alternatively, as demonstrated in \cite{Breiten2015,Zulfiqar2021}, using appropriately chosen weighting filters, the cost (\ref{eqn:cost}) can be approximated in a form that can be optimized in the computationally scalable interpolatory framework.   However, due to the approximations necessarily introduced in both classes of methods, there are no convergence guarantees with respect to the original problem.

Based on the above, it can be concluded that existing techniques are either not scalable to large-scale systems, or approximations introduced for scalability prevent convergence guarantees with respect to the original large-scale problem. Additionally, the methods often rely on the availability of state-space realizations of the dynamic systems - making their application to problems with non-realizable dynamics impossible. In this work, we aim to bridge this gap by proposing a novel method that is computationally scalable, has convergence guarantees and does not require state-space realizations. The method is inspired by analogous challenges encountered in the field of deep learning, where large-scale \emph{static} optimization problems are studied. In such a deep learning context, the vast amount of data and the increasingly deep neural networks involved also render exact computation of gradients intractable. An almost universally adopted solution to this is challenge in the field of deep learning is Stochastic Gradient Descent (SGD) \cite{bottou2012stochastic}. SGD substitutes exact gradients by stochastic estimates obtained by evaluating the cost function on randomly sampled subsets of the data. The practical success of SGD is supported by a large number of theoretical convergence guarantees, see \cite{Mertikopoulos2020,hardt2016train,vlaski2021second,wang2022uniform} for a small selection. We emphasize that the aforementioned methods are only applicable to static problems and that the problem in (\ref{eqn:cost}) involves the optimization of a dynamical system.

In this work, an extension of SGD for large-scale \emph{dynamical} $\mathcal{H}_2$ optimization problems is proposed. In the field of system theory, stochastic methods have been studied for the analysis of uncertain systems \cite{tempo2007monte,tempo2013randomized}. Although the problem setting in those works is different from ours, the choice for stochastic methods is similarly motivated by the recognition that viewing the problem in a stochastic setting can lead to computationally scalable methods. In \cite{peretz2016randomized}, a ray-shooting type stochastic method for finding minimal-gain static-output feedback controllers is proposed, and a specialization of the method for optimal PID controller design is presented in \cite{peretz2018randomized}. Probabilistic global convergence guarantees are provided - a second potential advantage of stochastic methods. Although the ray-shooting type methods demonstrate the potential of stochastic methods for optimization problems involving dynamical systems, they differ from the problem setting we consider in two aspects. First, we explicitly consider large-scale dynamical systems and develop a method that can efficiently handle sparse state-space representations. Second, even if no state-space representation of the dynamics are available, our method can be utilized as it relies only on frequency-domain samples.

The novel contributions presented in this paper are as follows:
\begin{enumerate}
	\item We propose an iterative stochastic method based on SGD to minimize $\hnorm{G(\mu)}$ for large-scale dynamical systems $G(\mu)$. In addition, probabilistic stability and convergence guarantees are provided for this approach.
	\item The method is easy to implement and requires only frequency-domain samples. Hence, it can be applied directly to problems for which no state-space realization is available.
	\item We show that the method is versatile and can be applied to a wide range of problems in the area of system theory by demonstrating it on two illustrative examples: 1) fixed-order observer design for a large-scale FEM thermal model, and 2) controller tuning for an infinite-dimensional plant.
\end{enumerate}

The remainder of this paper is organized as follows. First, the problem formulation is introduced in Section \ref{sct:problem_formulation}. A brief overview of SGD is presented in Section \ref{sct:sgd}. The proposed method is developed in Section \ref{sct:estimator_derivation} and its stability and convergence properties are analyzed in Section \ref{sct:convergence}. An algorithm implementing the method is given in Section \ref{sct:algorithm} along with an analysis of the computational complexity and practical suggestions for initialization and design of the stochastic sampling strategy. Finally, numerical examples are given in Section \ref{sct:numerical_examples} followed by concluding remarks in Section \ref{sct:concluding_remarks}.

We now discuss the notation used throughout the paper. The imaginary unit is denoted by $i$. $\stablereal$ denotes the set of linear, time-invariant, asymptotically stable and proper dynamical systems with rational transfer matrices. Note that the input-output dimension of elements of $\stablereal$ is not explicitly indicated. By a slight abuse of notation, we denote by $H(s)$ the transfer function associated with the dynamical system $H$ (and similarly, $G(\mu; s)$ is the transfer function of a parametric dynamical system $G(\mu)$).
$\left[ a \right]_i$ denotes the i-th component of a vector $a$. $\real\left(a\right)$ is the real part of the complex quantity $a$ and its complex conjugate is given by $\overline{a}$. The closed left-half of the complex plane is denoted by $\mathbb{C}_-$. Unless indicated otherwise, $\norm{ \cdot }$ represents the Euclidean norm for vector quantities.

\section{Problem formulation} \label{sct:problem_formulation}
In this work, we are interested in parameter optimization of the $\mathcal{H}_2$ system norm of a parameter-dependent continuous-time linear time-invariant (LTI) system $G(\mu)$, with $\mu \in \mathbb{R}^{n_{\mu}}$ the parameter. The optimization problem is formulated in terms of the cost function in (\ref{eqn:cost}),
where $\hnorm{\cdot}$ is induced by the $\mathcal{H}_2$ inner product defined as follows.
\begin{defn} \label{defn:h2inner}
The $\mathcal{H}_2$ inner product for two asymptotically stable and proper LTI systems $A$ and $B$ is
\begin{equation}
	\hinner{A}{B} = \frac{1}{2\pi} \int_{-\infty}^{\infty} \tr \left( A(i\omega) \tp{B(-i\omega)} \right) \mathrm{d}\omega.
	\label{eqn:h2inner}
\end{equation}
If one or both of $A$ and $B$ are unstable, then we set $\hinner{A}{B} = \infty$.
\end{defn}
The optimization problem associated with cost function (\ref{eqn:cost}) is
\begin{equation}
	c^{*} = \inf_{\mu \in \muSpace} c(\mu),
	\label{eqn:opt}
\end{equation}
where $\muSpace \subset \mathbb{R}^{n_{\mu}}$ is open and non-empty (implying existence of $c^{*}$). It is assumed that the number of degrees-of-freedom $n_{\mu}$ is low relative to the complexity of evaluating (\ref{eqn:cost}). This assumption is valid in many scenarios involving parameter optimization of large-scale dynamical systems, as highlighted by the examples in Fig. \ref{fig:optimization_examples} where $\mu$ parametrizes a low-dimensional entity such as a low-order model (a) or a reduced-order controller (b).

We explicitly do not assume that $G(\mu)$ is stable for all $\mu \in \muSpace$, i.e., the parameter-dependent system may be unstable for some choices of the parameter (note that in this case in accordance with Definition \ref{defn:h2inner} we set $c(\mu) = \infty$). We denote the set of values for $\mu$ for which $G(\mu)$ is asymptotically stable as $\muSpaceStab \subseteq \muSpace$ and adopt the following assumption that ensures a non-trivial solution to (\ref{eqn:opt}) exists.
\begin{assum} \label{assum:non_empty_muspacestab}
$\muSpaceStab$ is non-empty, i.e., there exists some $\mu \in \muSpace$ such that $G(\mu)$ is asymptotically stable.
\end{assum}

To obtain the stability and convergence guarantees of the proposed method that are derived in Section \ref{sct:convergence}, we additionally adopt the following assumptions on the regularity of $c(\mu)$.
\begin{assum} \label{assum:c_continuity}
$c$ is continuous for all $\mu \in \muSpace$.
\end{assum}
\begin{assum}
$c$ restricted to any sub-level set is Lipschitz continuous and Lipschitz smooth, i.e., for each $C < \infty$, we have that there exist $K, L < \infty$ such that
$$
| c(\mu) - c(\mu') | \leq K \norm{\mu - \mu'},
$$
and
$$
\norm{ \nabla c(\mu) - \nabla c(\mu') } \leq L \norm{\mu - \mu'},
$$
for all $\mu, \mu' \in \{ \mu | c(\mu) \leq C \}$.
\label{assum:c_sublevel_lipschitz}
\end{assum}
\begin{assum}
$c$ is coercive, i.e., $c(\mu) \rightarrow \infty$ as $\norm{\mu} \rightarrow \infty$.
\label{assum:c_coercivity}
\end{assum}
Informally, Assumption \ref{assum:c_continuity} ensures that we ``do not go to instability without noticing". It means that the cost function of a trajectory of $\mu$ transitioning from a stable parameter to an unstable parameter always goes to $\infty$ as the unstable parameter region is approached. This does not necessarily hold for arbitrary $G(\mu)$ due to unstable pole cancellations, as can be observed from the simple example
\begin{equation}
	G(\mu; s) = \frac{2 \mu}{s + \mu}
\end{equation}
with $\muSpace = \mathbb{R}$ for which the cost function is
\begin{equation}
	c(\mu) = \begin{cases}
		\mu, \text{ if } \mu \geq 0, \\
		\infty, \text{ if } \mu < 0.
	\end{cases}
\end{equation}
In Appendix \ref{sct:continuity}, we characterize a class of parameter-dependent systems for which Assumption \ref{assum:c_continuity} is indeed guaranteed to hold. Assumption \ref{assum:c_sublevel_lipschitz} is another regularity requirement that enables the guarantee on stability preservation of the proposed method which is derived later in the paper. Finally, Assumption \ref{assum:c_coercivity} together with Assumptions \ref{assum:c_continuity} and \ref{assum:c_sublevel_lipschitz} ensures that sub-level sets of $c(\mu)$ are compact and contain at least one minimizer. It may fail to hold in scenarios such as optimal control problems (Figure \ref{fig:optimization_examples:synthesis}) where actuation effort is not penalized. As a result, actuation effort can be indefinitely increased ($\norm{\mu} \rightarrow \infty$, where $\mu$ could reflect controller gains) to improve performance. Thus, to satisfy Assumption \ref{assum:c_coercivity} it is important all relevant input-output relations of the underlying system-to-be-optimized appear in the cost. Alternatively, one may work with a modified cost function $c'(\mu) := c(\mu) + c_\mathrm{reg}(\mu)$ involving a regularization term $c_\mathrm{reg}$ which has the desired property $c_\mathrm{reg}(\mu) \rightarrow \infty$ as $\norm{\mu} \rightarrow \infty$. 

Assumption \ref{assum:c_sublevel_lipschitz} implies existence of the gradient $\nabla c(\mu)$ for all $\mu \in \muSpaceStab$. This gradient can be expressed in terms of the $\mathcal{H}_2$ inner product (\ref{eqn:h2inner}) as follows:
\begin{equation}
	\nabla c(\mu) = \tp{\begin{bmatrix} {\hinner{G(\mu)}{\dGdmu{1}(\mu)}} & \cdots & {\hinner{G(\mu)}{\dGdmu{n_\mu}(\mu)}} \end{bmatrix}}.
	\label{eqn:nabla_c}
\end{equation}
If state-space descriptions of $G$ and its parameter-gradient systems are available, then the exact gradient can be found \cite[Lemma 2.1.5]{Antoulas2020}. However, it involves computing the controllability or observability Gramians of $G$ and the associated parameter-gradient systems, which scales with $n^3$ (where $n$ is sum of the state-space dimensions of $G$ and its parameter-gradient systems). Hence, this is computationally infeasible for the large-scale systems (e.g., with $n$ in the order of millions) considered in this work. Thus, an approximation of the gradient is required. However, the error introduced by this approximation will typically also lead to an approximate solution of the optimization problem in (\ref{eqn:opt}) with no guarantees on the quality of the approximate solution. Inspired by results of stochastic gradient descent (SGD) in the machine learning community, we are interested in approximating the gradient in such a way that we can still give probabilistic guarantees on the quality of the solution to (\ref{eqn:opt}). This leads to the following problem formulation:

\emph{How can we efficiently approximate the gradient (\ref{eqn:nabla_c}) such that we can still guarantee (probabilistic) convergence to a local minimizer of (\ref{eqn:opt})?}

\section{Stochastic gradient descent in static problems} \label{sct:sgd}

In this section, we recall how stochastic gradient descent (SGD) renders large-scale \emph{static} optimization problems feasible. We use the term static here to emphasize the fact that the optimization does not involve dynamical system models. As a stepping stone for the novel results in Section \ref{sct:convergence}, we briefly discuss historical developments of SGD and describe some of the convergence guarantees available in the literature. For a more complete overview consult, e.g., \cite{Netrapalli2019}.

The origins of SGD can be traced back to the more general framework of stochastic approximation techniques, which involve the iterative optimization of functions which cannot be evaluated directly (for example, because only noisy data is available or an exact evaluation is computationally infeasible), but for which some random estimator is given. The first stochastic approximation technique is the Robbins-Monro algorithm \cite{Robbins1951}. It is intended as an iterative method for finding the roots of a function $M(x)$ which cannot be evaluated exactly (for example, because it is computationally infeasible). However, it is assumed that random estimates $N(x)$ of $M(x)$ are available, with the property that the measurements are unbiased, i.e., $\mathbb{E}[N(x)] = M(x)$. A link between the Robbins-Monro algorithm and gradient-based minimization of a differentiable function $f(x)$ can be established by observing that local minimizers are roots of the equation $\nabla f(x) = 0$. Solving this problem using random estimates of $\nabla f(x)$ is the core concept of SGD. One area in which SGD has exploded in popularity in recent years is machine learning \cite{Bottou2010,Ge2021}. This popularity may be explained in two parts. First, the models and data involved in machine learning problems have grown exponentially in size, prohibiting exact calculation of gradients. Second, in many applications it is straightforward to design a random estimator of the gradient whose statistical properties provide attractive convergence guarantees. The latter point will be explained in more detail next.

In machine learning problems, the cost function $c(\mu)$ is often a sum of the contributions of each element of the dataset $c_{\mathrm{element}}$:
\begin{equation}
	c(\mu) = \sum_{i=1}^{N} c_{\mathrm{element}}(\mu; x_i),
\end{equation}
with $N$ the dataset size, $x_i$ element $i$ of the dataset and $\mu$ the parameters of the model (e.g., neural network) to be optimized. The gradient is also a sum of contributions:
\begin{equation}
	\nabla c(\mu) = \sum_{i=1}^{N} \nabla c_{\mathrm{element}}(\mu; x_i).
\end{equation}
Since $N$ is usually large, calculating the exact gradient is computationally too expensive. Instead, the (batch) stochastic gradient is calculated by selecting at each iteration a small random subset $\mathcal{I}$ of the dataset elements, and using that as an estimate of the gradient:
\begin{equation}
	\widehat{\nabla c}(\mu) = \sum_{i \in \mathcal{I}} \nabla c_{\mathrm{element}}(\mu; x_i).
	\label{eqn:nabla_c_estimator_static}
\end{equation}
The gradient estimate (\ref{eqn:nabla_c_estimator_static}) replaces the true gradient in whichever procedure is selected to optimize the parameters $\mu$. For example, in this paper, we consider gradient descent:
\begin{equation}
	\mu_{k+1} = \mu_{k} - \alpha_k \widehat{\nabla c}(\mu_k)
	\label{eqn:sgd_static}
\end{equation}
with $k$ the iteration counter and $\alpha_k$ the step size. Repeated application of (\ref{eqn:sgd_static}) until some pre-specified convergence criterion is achieved, forms the basic SGD algoritm. Although this form of SGD is simple, it is effective in practice \cite{Bottou2010}. Furthermore, its convergence properties for both convex and non-convex optimization problems has been studied in literature \cite{Mertikopoulos2020,Lei2019,So2017}. Generally speaking, it can be shown that using (\ref{eqn:nabla_c_estimator_static}) as a surrogate for the exact gradient in (\ref{eqn:sgd_static}) in the static setting has attractive convergence guarantees because the gradient error $Z(\mu) := \nabla c(\mu) - \widehat{\nabla c}(\mu) $ satisfies two key stochastic properties (recall that the subset $\mathcal{I}$ of dataset elements is sampled randomly at each iteration of SGD):
\begin{enumerate}
	\item[(P1)] \textbf{Unbiased}: $\mathbb{E}\left[ Z(\mu) \right] = 0$,
	\item[(P2)] \textbf{Finite variance}: for some $\sigma < \infty$ we have $\mathbb{E}\left[ \norm{Z(\mu)}^2 \right] \leq \sigma^2, \quad \forall \mu \in \muSpace$.
\end{enumerate}
Although the practical success of SGD greatly depends on the ability to derive estimators that are both unbiased and have low variance for a given problem, several other types of improvements to SGD have been proposed, such as leveraging second-order information to speed up convergence \cite{Bordes2009} and implicit SGD updates to decrease sensitivity to hyperparameters \cite{Toulis2017}.

As discussed in this section, in the context of static optimization problems where the amount of data and complexity of the models is growing rapidly, SGD is an important ingredient to ensure the optimization remains computationally feasible while at the same time maintaining convergence guarantees. However, SGD approaches for dynamical problems are not available in literature. Motivated by its advantages, we explore extending SGD to problems involving parametric dynamical systems in the next section.

\section{Extending SGD to dynamical problems} \label{sct:estimator_derivation}
The fundamental ingredient of any SGD approach is an effective random estimator of the gradient (such as the estimator (\ref{eqn:nabla_c_estimator_static}) in the context of static problems). Special care needs to be taken that the stochastic properties discussed above (unbiased and finite variance) hold, since they are crucial to derive convergence guarantees. In this section, a gradient estimator for the dynamical optimization problem (\ref{eqn:opt}), which is both computationally inexpensive and meets these stochastic properties, will be proposed.

As the first step, using (\ref{eqn:h2inner}) and (\ref{eqn:nabla_c}) the gradient is written in a frequency-domain integral form:
\begin{equation}
	\nabla c(\mu) = \int_{-\infty}^{\infty} f(\mu; i\omega) d\omega.
	\label{eqn:nabla_c_j}
\end{equation}
where $f(\mu)$ represents a single-input multiple-output parameter-dependent dynamical system with transfer function $f(\mu; s)$ for which the j-th output is 
\begin{equation}
	[f(\mu)]_j = \frac{1}{2\pi} \mathrm{tr}\left( G(\mu) \overline{\frac{\partial G}{\partial \mu_j}}(\mu)^{\mathrm{T}} \right).
	\label{eqn:f}
\end{equation}
As mentioned previously, evaluation of the exact gradient is computationally intractable for large-scale systems. However, the integral form in (\ref{eqn:nabla_c_j}) allows us to employ the technique of \emph{Monte Carlo integration} \cite{elvira2014advances}. In general, Monte Carlo integration approximates an integral by sampling the integrand at a discrete and finite set of points that are drawn from some probability distribution. Subsequently, a weighted mean of the integrand at these points delivers the approximation. Let frequency points $\{\omega_m\}_{m=1}^{M}$ be $M$ independent random variables drawn from a Probability Density Function (PDF) $p(\omega)$ with $\omega \in \mathbb{R}$. Using Monte Carlo integration an estimator of (\ref{eqn:nabla_c_j}) is given by
\begin{equation}
	\widehat{\nabla c}(\mu) = \frac{1}{M} \sum_{m=1}^{M} \frac{f(\mu; i\omega_m))}{p(\omega_m)}.
	\label{eqn:nabla_c_j_estimator}
\end{equation}
The function $f(\mu, \cdot)$ in (\ref{eqn:nabla_c_j}) and (\ref{eqn:nabla_c_j_estimator}) is complex-valued. However, the true integral is  real-valued since $f(\mu)$ represents a real dynamical system, implying conjugate symmetry of $f$ ($f(\mu; i\omega) = \overline{f(\mu; -i\omega)}$). Thus, to ensure our estimate is always real as well, we exploit this conjugate symmetry by choosing $p(\omega)$ such that $p(\omega) = p(-\omega)$ and obtain
\begin{equation}
	\widehat{\nabla c}(\mu) = \frac{1}{M} \sum_{m=1}^{M} \frac{\real \left( f(\mu; i\omega_m) \right)}{p(\omega_m)}.
	\label{eqn:nabla_c_hat}
\end{equation}

In Sections \ref{sct:p1}-\ref{sct:p2}, we analyse under which conditions the estimator (\ref{eqn:nabla_c_hat}) satisfies requirements (P1) and (P2).

\subsection{(P1) Unbiased estimate} \label{sct:p1}
Let $\Omega := \{ \omega \in \mathbb{R} | p(\omega) > 0 \}$ and denote by $\Omega^{\mathrm{c}}$ its complement in $\mathbb{R}$. We aim to find the expected value of the estimator (\ref{eqn:nabla_c_hat}) to determine its bias with respect to the true gradient. The estimator is a function of the random variables $\omega_m$, $m = 1, ..., M$, which are independently drawn from the same distribution (namely, the PDF given by $p(\omega)$). Thus, the expected value of each term of the sum in the right-hand side of (\ref{eqn:nabla_c_hat}) is equal and given by
\begin{equation}
	\mathbb{E}\left[ \frac{ \mathrm{Re}(f(\mu; i\omega_m)) }{p(\omega_m)} \right] = \int_{\Omega} \frac{ \mathrm{Re}(f(\mu; \omega)}{p(\omega)} p(\omega) d\omega,
	\label{eqn:estimator_expected_value}
\end{equation}
for all $m = 1, ..., M$. Making use of this fact, we have
\begin{equation}
	\begin{aligned}
	& \mathbb{E}\left[ \widehat{\nabla c}(\mu) \right]
 = \mathbb{E}\left[ \frac{1}{M} \sum_{m=1}^{M} \frac{\real \left( f(\mu; i\omega_m) \right)}{p(\omega_m)} \right] \\
 & = \frac{1}{M} \sum_{m=1}^{M} \mathbb{E}\left[ \frac{\real \left( f(\mu; i\omega_m) \right)}{p(\omega_m)} \right] \\
	& = \frac{1}{M} \sum_{m=1}^{M} \int_{\Omega} \frac{\real \left( f(\mu; i\omega) \right)}{p(\omega)} p(\omega) d\omega \\
	& = \int_{\Omega} \real \left( f(\mu; i\omega) \right)  d\omega \\
	& = \int_{-\infty}^{\infty} \real \left( f(\mu; i\omega) \right)  d\omega - \int_{\Omega^{\mathrm{c}}} \real \left( f(\mu; i\omega) \right)  d\omega \\
	& = \nabla c(\mu) - \int_{\Omega^{\mathrm{c}}} \real \left( f(\mu; i\omega) \right)  d\omega.
	\end{aligned}
	\label{eqn:chat_mean}
\end{equation}
In (\ref{eqn:chat_mean}), the third equality follows from the fact that the random variables $\omega_m$, $m = 1, ..., M$, are independent and identically distributed. The bias of the gradient estimate is then given by
\begin{equation}
	\mathbb{E}[ Z(\mu) ] = \int_{\Omega^{\mathrm{c}}} \real \left( f(\mu; i\omega) \right) d\omega.
	\label{eqn:bias}
\end{equation}
Compared to SGD for static problems, the estimator (\ref{eqn:nabla_c_hat}) is unbiased only if the following conditions are satisfied:
\begin{enumerate}
	\item $\mu \in \muSpaceStab$, that is, $G(\mu)$ is asymptotically stable. Now, the following issue arises: even though we have not defined the gradient for $\mu \notin \muSpaceStab$, frequency samples of $f(\mu; \cdot)$ remain finite with probability 1 (since poles of $f(\mu; \cdot)$ are sampled with probability 0) and the estimator is thus unable to detect instability.
	\item The integral (\ref{eqn:bias}) is zero. To guarantee this, frequency regions where $\mathrm{Re}(f(\mu; \cdot))$ is non-zero should have a non-zero probability of being sampled.
\end{enumerate}
The first condition highlights the importance of preserving stability of the underlying system during the iterations of SGD - because only then can we expect meaningful convergence guarantees. We will discuss stability preservation in detail in Section \ref{sct:convergence}.

\subsection{(P2) Bounded variance} \label{sct:p2}
Property (P2) is satisfied if and only if the error of the gradient estimate has bounded variance, i.e., if and only if there exists a $\sigma < \infty$ such that
\begin{equation}
	\mathbb{E}\left[ \norm{\nabla c(\mu) - \widehat{\nabla c}(\mu)}^2 \right] < \sigma^2, \quad \forall \mu \in \muSpaceStab,
\end{equation}
where we are, similar to the analysis on bias in the previous section, again limited to the asymptotically stable part of the parameter space $\muSpace$. However, if $\muSpace \setminus \muSpaceStab$ is non-empty (i.e., there exist $\mu$ making $G(\mu)$ unstable), then $\nabla c(\mu)$ is unbounded on $\muSpaceStab$ due to Assumption \ref{assum:c_continuity}. As a result, it is in general not reasonable to expect that the variance of the error is bounded.

To circumvent this issue, we first proceed by showing (in Theorem \ref{thm:bounded_variance}) under which condition on the PDF $p(\omega)$ we can guarantee bounded variance on an arbitrary sub-level set of $c(\mu)$. However, in order to then prove convergence of the SGD iterates (\ref{eqn:sgd_static}) later in the paper, we need the iterands $\mu_n$, $n = 1, 2, ...$, to remain in this sub-level set (so that we have bounded variance for all iterands). Thus, in Theorem \ref{thm:stability} in Section \ref{sct:convergence}, we prove that the iterands $\mu_n$, $n = 1, 2, ...$, of SGD indeed remain within this sub-level set with high probability.
\begin{thm} \label{thm:bounded_variance}
For any $C < \infty$, there exists a $\sigma < \infty$ such that
$$
\mathbb{E}\left[ \norm{\nabla c(\mu) - \widehat{\nabla c}(\mu)}^2 \right] < \sigma^2, \quad \forall \mu \in \{ \mu | c(\mu) \leq C \},
$$
with $\nabla c(\mu)$ as in (\ref{eqn:nabla_c}) and $\widehat{\nabla c}(\mu)$ as in (\ref{eqn:nabla_c_hat}), if the sampling PDF $p$ used in (\ref{eqn:nabla_c_hat}) satisfies for all $\omega \in \Omega$
\begin{equation}
\frac{\norm{\real\left(f(\mu; i\omega)\right)}}{\sqrt{p(\omega)}} = \mathcal{O}\left(\frac{1}{|\omega| + 1}\right).
\label{eqn:pdf_bound}
\end{equation}
\end{thm}

\begin{pf}
Using the triangle inequality, we have
\begin{equation}
	\begin{aligned}
		& \mathbb{E}\left[ \norm{ \nabla c(\mu) - \widehat{\nabla c}(\mu)}^2 \right] \\
		& \leq \mathbb{E}\left[ \norm{ \widehat{\nabla c}(\mu) }^2 \right] + \norm{ \nabla c(\mu) }^2
	\end{aligned}
	\label{eqn:variance_step1}
\end{equation}
Note that $\norm{ \nabla c(\mu) }$ is bounded (due to Assumption \ref{assum:c_sublevel_lipschitz}). Thus, it is sufficient to prove that the first term on the right-hand side of (\ref{eqn:variance_step1}) is bounded. It can be written as follows:
\begin{equation}
	\mathbb{E}\left[ \norm{ \widehat{\nabla c}(\mu) }^2 \right] = \mathbb{E}\left[ \norm{ \frac{1}{M} \sum_{m=1}^{M} \frac{\real(f(\mu; i\omega_m))}{p(\omega_m)} }^2 \right].
	\label{eqn:variance_step2}
\end{equation}
Introduce the short-hand notation
$$A_m := \frac{\mathrm{Re}(f(\mu; i\omega_m))}{p(\omega_m)}, \quad m = 1, ..., M.$$ The random variables $A_m$ are independent and identically distributed (see also (\ref{eqn:estimator_expected_value})), i.e., $\mathbb{E}[ A_i^{T} A_j ] = \mathbb{E}[ A_i ]^{T} \mathbb{E}[ A_j ] = \norm{\mathbb{E}[A_1]}^2$. Making use of this and expanding (\ref{eqn:variance_step2}) using this short-hand notation gives
\begin{equation}
\begin{aligned}
	& \mathbb{E}\left[ \norm{ \widehat{\nabla c}(\mu) }^2 \right] = \mathbb{E}\left[ \norm{ \frac{1}{M} \sum_{m=1}^{M} A_m }^2 \right] \\
	& = \frac{1}{M^2} \left( M \mathbb{E}\left[ \norm{ A_1 }^2 \right] + M(M-1) \norm{ \mathbb{E}[ A_1 ]}^2 \right) \\
	& = \frac{1}{M} \left( \mathbb{E}\left[ \norm{ A_1 }^2 \right] + (M-1) \norm{ \mathbb{E}[ A_1 ]}^2 \right).
\end{aligned}
	\label{eqn:variance_step3}
\end{equation}
Neglecting the factors depending on $M$ (which are finite), we obtain for the first term on the right-hand side of (\ref{eqn:variance_step3}) the following:
\begin{equation}
\begin{aligned}
	\mathbb{E}\left[ \norm{A_1}^2 \right] & = \int_{\Omega} \frac{\norm{ \mathrm{Re}(f(\mu; i\omega))}^2}{p(\omega)^2} p(\omega) d\omega \\
	& = \int_{\Omega} \frac{\norm{ \mathrm{Re}(f(\mu; i\omega))}^2}{p(\omega)} d\omega,
\end{aligned}
\end{equation}
which is bounded if the integrand is $\mathcal{L}_2$ integrable, which, in turn, is guaranteed if (\ref{eqn:pdf_bound}) holds. For the second term of (\ref{eqn:variance_step3}) we have
\begin{equation}
\begin{aligned}
	\norm{\mathbb{E}\left[A_1\right]} & = \norm{ \int_{\Omega} \frac{ \mathrm{Re}(f(\mu; i\omega))}{p(\omega)} p(\omega) d\omega } \\
	& = \norm{ \int_{\Omega} \mathrm{Re}(f(\mu; i\omega)) d\omega } \\
	& \leq \int_{\Omega} \norm{ f(\mu; i\omega)} d\omega \\
	& \leq \int_{-\infty}^{\infty} \norm{ f(\mu; i\omega)} d\omega,
\end{aligned}
\end{equation}
which is bounded since the stability properties of $f(\mu)$ are inherited from $G(\mu)$ (see (\ref{eqn:f})) and $\hnorm{G(\mu)}$ is bounded (since $\mu \in \{ \mu | c(\mu) \leq C \}$).
 $\blacksquare$
\end{pf}
Theorem (\ref{thm:bounded_variance}) provides a sufficient condition on the PDF $p(\omega)$ used in the gradient estimator (\ref{eqn:nabla_c_hat}) such that we have bounded variance of the gradient estimate error on sub-level sets of $c(\mu)$. In Section \ref{sct:pj}, we give some examples to show how to construct PDFs satisfying this condition.

\section{Convergence guarantees} \label{sct:convergence}
The previous analyses indicated that the proposed estimator (\ref{eqn:nabla_c_hat}) is unbiased only if $\mu \in \muSpaceStab$. Furthermore, bounded variance is only guaranteed if $\mu$ remains within a fixed sub-level set of $c(\mu)$. Since unbiasedness and bounded variance are desirable properties in the scope of achieving convergence results, it is important to study conditions under which we can guarantee that the cost $c(\mu_k)$ remains bounded for all SGD iterands $\mu_k$, $k = 0, 1, 2, ...$, i.e., for which we can guarantee that these iterands all relate to asymptotically stable systems. In literature, several boundedness results have been derived. In \cite{wang2022uniform}, uniform boundedness of iterates in expectation is guaranteed for a class of non-convex functions with sufficiently fast asymptotic growth. However, guarantees in expectation are less useful here because as soon as an iterate enters the unstable regime, the estimator becomes biased and therefore convergence guarantees on subsequent iterates are no longer applicable. Instead, we are more interested in showing that the probability that \emph{all} iterates remain stable can be controlled. In fact, the following theorem is the main theoretical contribution of this work and shows that this probability can be moved arbitrarily close to 1 by an appropriate step size policy (i.e., a policy for $\alpha_k$ in the SGD scheme in (\ref{eqn:sgd_static})).
\begin{thm}
Let $G(\mu)$ be a dynamical system for which $c(\mu) = \hnorm{G(\mu)}^2$ satisfies Assumptions \ref{assum:non_empty_muspacestab}, \ref{assum:c_continuity} and \ref{assum:c_sublevel_lipschitz}. Given an initial parameter $\mu_0 \in \muSpaceStab$ ($c(\mu_0) < \infty$) and tolerance $\varepsilon > 0$, define the event that all iterates remain bounded:
\begin{equation}
E_\infty := \{ c(\mu_k) \leq c(\mu_0) + \sqrt{\varepsilon} + \varepsilon | k = 1, 2, ... \}.
\label{eqn:event_stable}
\end{equation}
Then, for any tolerance $\delta > 0$, we have that
\begin{equation}
	\mathbb{P}[E_\infty] \geq 1 - \delta
	\label{eqn:Einfty_bound}
\end{equation}
if the step size policy $\alpha_k$ satisfies
\begin{equation}
	\sum_{k=0}^{\infty} \alpha_k^2 \leq \frac{\delta \cdot \varepsilon}{R_*}
	\label{eqn:alpha_sum_bound}
\end{equation}
and $0 \leq \alpha_k \leq L^{-1}$ for all $k = 0, 1, ...$, with
\begin{equation}
	R_* = K^2 \sigma^2 + 2L(K^2 + \sigma^2)
	\label{eqn:Rstar}
\end{equation}
and $K, L$ the Lipschitz constants (Assumption \ref{assum:c_sublevel_lipschitz}) and $\sigma$ an upper bound on the variance (Theorem \ref{thm:bounded_variance}) on the sublevel set $\{ \mu | c(\mu) \leq c(\mu_0) + \sqrt{\varepsilon} + \varepsilon \}$, respectively.
\label{thm:stability}
\end{thm}
\begin{pf}
See Appendix \ref{sct:boundedness_of_iterates}.
\end{pf}
To the authors' knowledge, this is the first result guaranteeing with high probability boundedness of the $\mathcal{H}_2$ norm of a dynamical systems $G(\mu)$ (and thus the asymptotic stability of the dynamical system) under SGD iterates. Importantly, it provides a solid theoretical foundation for applying SGD to problems involving dynamical systems. Namely, conditioned on the event (\ref{eqn:event_stable}) that the iterands stay bounded, many convergence guarantees from literature for static problems can be applied to our proposed SGD scheme for dynamical problems. In the remainder of this section, we will provide 2 examples of such results. First, the following theorem shows that the proposed SGD scheme converges almost surely to a stationary point.
\begin{thm} \label{thm:convergence_to_stationary_point}
If Assumption \ref{assum:c_coercivity} holds and the step size satisfies $\alpha_k = \Theta(1/k^p)$ for some $p \in (1/2, 1]$, then conditioned on the event $E_\infty$ as in (\ref{eqn:event_stable}), the sequence $\mu_k$, $k=1, 2, ...$ almost surely converges to a point $\mu^*$ such that $\nabla c(\mu^*) = 0$.
\end{thm}
\begin{pf}
Since we are conditioned on $E_\infty$, we know that $c(\mu_k)$ remains within a sub-level set. Assumption \ref{assum:c_coercivity} implies that this sub-level set is bounded. Hence, we have satisfied the conditions for applying \cite[Theorem 2]{Mertikopoulos2020}. $\blacksquare$
\end{pf}
Note that the step size policy $\alpha_k = \Theta(1/k^p)$, $p \in (1/2, 1]$, in Theorem \ref{thm:convergence_to_stationary_point} is $\ell_2$ summable and, hence, is compatible with the step size restriction (\ref{eqn:alpha_sum_bound}) for stability preservation.

Stationary points ($\nabla c(\mu) = 0$) may correspond to saddle points, which we would like to avoid when seeking for a (local) minimizer. Whereas deterministic gradient descent is unable to escape from a stationary point, the stochastic perturbations present in the SGD updates can help us avoid stationary points that do not correspond to local minimizers. Strong guarantees can be given in the case of \emph{strict} saddle points (i.e., points where the minimum eigenvalue of the Hessian is smaller than 0).
We will formalize this by first introducing the concept of strict saddle manifolds \cite{Mertikopoulos2020}.
\begin{defn} \label{defn:strict_saddle_manifold}
A \emph{strict saddle manifold} $\mathcal{S} \subset \muSpace$ of $c(\mu)$ is a smooth connected component of the set of stationary points $\{ \mu \in \muSpace | \nabla c(\mu) = 0 \}$ such that the following two properties hold:
\begin{enumerate}
	\item Every $\mu \in \mathcal{S}$ is a strict saddle point of $c(\mu)$.
	\item There exists a $\tau > 0$ such that for all $\mu \in \mathcal{S}$ all non-zero eigenvalues $\lambda$ of the Hessian satisfy $| \lambda | \geq \tau$.
\end{enumerate}
\end{defn}
We can then show that the proposed method avoids strict saddle manifolds with probability 1 with the following theorem.
\begin{thm}
Let the step size $\alpha_k = \Theta(1/k^p)$ for some $p \in (0, 1]$. Then, conditioned on the event $E_\infty$ in (\ref{eqn:event_stable}), we have that
\begin{equation}
	\mathbb{P}\left[ \lim_{k \rightarrow \infty} \mu_k \in \mathcal{S} \right] = 0,
\end{equation}
with $\mathcal{S}$ a strict saddle manifold, if there exists a $\gamma > 0$ such that for all $\mu$ contained in the sub-level set associated with $E_\infty$ in (\ref{eqn:event_stable}) and all $b \in \mathbb{R}^{n_\mu}$ with $\norm{b}_2 = 1$ we have
\begin{equation}
	\mathbb{E}\left[ | Z(\mu)^\mathrm{T} b | \right] \geq \gamma.
	\label{eqn:Z_uniform_exciting}
\end{equation}
\end{thm}
\begin{pf}
See \cite[Theorem 3]{Mertikopoulos2020}.
\end{pf}
The condition (\ref{eqn:Z_uniform_exciting}) states that the gradient error $Z(\mu)$ has non-zero component in every direction and emphasizes that the result does not hold if some or all of the gradient components are deterministic estimates. The uniform bound $\gamma$ is typically automatically satisfied in practice or can be ensured by adding small artificial stochastic perturbations sampled from the unit sphere \cite{ge2015escaping}.

\section{Algorithmic implementation} \label{sct:algorithm}
An implementation of the SGD method for the stochastic $\mathcal{H}_2$ optimization of large-scale parametric dynamical systems (\shopt) proposed in this work is given in Algorithm \ref{alg:sh2opt}. Connected to this algorithm, in this section we discuss the practical issues of initialization $\mu_0$ (Section \ref{sct:initialisation}), choosing the sampling PDF $p(\omega)$ (Section \ref{sct:pj}) and computational complexity (Section \ref{sct:computational_cost}).

\begin{algorithm}
	\caption{\shopt: Stochastic $\mathcal{H}_2$ optimization of dynamical systems}
	\label{alg:sh2opt}
	\begin{algorithmic}[1] 
		\State Given: parametrized dynamical system $G(\mu)$, parameter space $\muSpace$, initialization $\mu_0 \in \muSpaceStab$, iteration limit $N$, step sizes $\alpha_k \leq L^{-1}$ satisfying (\ref{eqn:alpha_sum_bound}), PDF $p(\omega)$, number of samples $M$
		\State Initialize iteration counter $k = 0$
		\LineComment{Loop over iterations}
		\While{$k \leq N$}
			\LineComment{Monte Carlo integration}
			\State Draw $M$ samples $\{ \omega_m \}_{m=1}^{M}$ from $p(\omega)$
			\State Evaluate $f(\mu_k; i\omega_m)$, $m = 1, ..., M$, using (\ref{eqn:f})
			\State $\widehat{\nabla c}(\mu_k) = \sum_{m=1}^{M} \frac{\mathrm{Re}(f\left(\mu_k; i\omega_m\right))}{p(\omega_m)}$
			\LineComment{Descent step}
			\State $\mu_{k+1} = \mu_k - \alpha_k \widehat{\nabla c}(\mu_k)$
			\State $k \leftarrow k + 1$
		\EndWhile
	\end{algorithmic}
\end{algorithm}

\subsection{Initialization} \label{sct:initialisation}
A good initialization $\mu_0$ for Algorithm \ref{alg:sh2opt} satisfies two properties:
\begin{enumerate}
	\item[(A)] $G(\mu_0)$ is asymptotically stable,
	\item[(B)] $\mu_0$ is close to a (local) minimizer of (\ref{eqn:cost}).
\end{enumerate} 
Regarding property (A), in some problems stability is guaranteed by the parametrization. For instance, if a passive plant is controlled through negative feedback by a passive controller, the closed-loop system is guaranteed to be passive (and hence stable). A parametrization $K(\mu)$ of passive controllers would therefore be a good option for such a problem. In cases where such guarantees are unavailable, an estimation of $\hnorm{G(\mu_0)}$ through the use of sparse Lyapunov equation solvers can be used \cite{Kurschner2016}. (An approximation of the cost function obtained with such a solver can also be used to terminate Algorithm \ref{alg:sh2opt}.)

Regarding property (B), an initialization can be found using a typical ``reduce-then-design" approach that is common in indirect fixed-order control design schemes. In this approach, the interconnected components forming $G(\mu)$ are first disconnected. In the case of control design (Figure \ref{fig:optimization_examples:synthesis}), this gives the plant model. Then, the open-loop components are reduced using standard $\mathcal{H}_2$ optimal approaches such as IRKA \cite{Gugercin2008} to the desired order. Finally, the reduced components are used to obtain the initialization $\mu_0$. In the case of control design, an LQG controller synthesized on the reduced plant model could be used for such purpose.

\subsection{Designing the sampling distribution $p$} \label{sct:pj}
The question how one should design the sampling distribution $p$ for the SGD can be steered by reviewing the various results that were derived in this article:
\begin{itemize}
	\item The estimator should be unbiased (for $\mu \in \muSpaceStab$).
	\item The estimator should have low variance (for $\mu \in \muSpaceStab$) which facilitates a larger step size (due to Theorem \ref{thm:stability}, see (\ref{eqn:alpha_sum_bound}), (\ref{eqn:Rstar})) and faster convergence.
\end{itemize}
A good starting point is to reflect on whether it is possible to have a sampling distribution that achieves both of these conditions perfectly (for a fixed $\mu$), i.e., is unbiased with \emph{zero} variance. In fact, this is possible if $\mathrm{Re}(f(\mu; i\omega))$ is scalar, does not change sign and satisfies
\begin{equation}
	|\mathrm{Re}(f(\mu; i\omega))| = \mathcal{O}\left( \frac{1}{|\omega|^2 + 1} \right).
	\label{eqn:f_condition}
\end{equation}
Then, the optimal sampling distribution $p^*$ achieving zero variance is $p^*(\omega) = \tau |\mathrm{Re}(f(\mu; i\omega))|$ where $\tau$ is a normalizing factor ensuring the distribution is valid \cite{Kahn1953}. The condition (\ref{eqn:f_condition}) ensures that $p^{*}$ satisfies (\ref{eqn:pdf_bound}).
Moreover, in the case that $\mathrm{Re}(f(\mu; i\omega)$ does change sign, $p^*$  still provides the lowest variance among all possible distributions \cite[Section 9.1]{Owen2013}.

One immediate issue preventing directly using $p^*$ is that the exact gradient is required (namely, to determine $\tau$). Hence, any computational benefit of the proposed SGD method would be lost. However, we can use $p^*$ to derive a rule-of-thumb for sampling distributions that are computationally attractive and provide low variance: the distribution should be approximately proportional to $|\mathrm{Re}(f(\mu; i\omega))|$. Prior knowledge on $f(\mu)$ can be exploited to achieve this. For example, if it is known that $G(\mu)$ and its parameter-gradient systems are band-limited in a certain frequency range, then the sampling distribution could be taken by constructing a low-order band-pass filter $B$ around this frequency range and then taking $p(\omega) = \tau |\mathrm{Re}(B(i\omega))|$. Alternatively, one could take several frequency samples of $G(\mu)$ and its parameter-gradient systems and let $B$ be an interpolant of these frequency samples (using, e.g., the AAA algorithm \cite{Nakatsukasa2018} in the single-input single-output case).

Simpler distributions are possible and can even be band-limited (at the expense of introducing bias). For example, if the system-to-be-optimized is (approximately) band-limited one may also sample uniformly or log-uniformly in that band. Drawing frequency samples from such distributions can be done very effectively as described in Section \ref{sct:computational_cost}, but the variance may be higher than when an interpolant-based distribution is used because the latter better captures the underlying rational structure of the transfer function $f(\mu; s)$. Another immediate benefit of utilizing a band-limited PDF is that condition (\ref{eqn:pdf_bound}) is trivially satisfied since $\Omega$ is bounded.

\subsection{Computational cost} \label{sct:computational_cost}
Recall that exact computation of the gradient (\ref{eqn:nabla_c}) is expensive for large-scale dynamical systems. It involves computation of the observability or controllabilty Gramians, requiring $\mathcal{O}(n^3)$ computational time (with $n$ being the total dimension of the state-space representations of $G(\mu)$ and $\frac{\partial G}{\partial \mu_j}(\mu)$, $j = 1, ..., n_\mu$). The computational cost of the estimator proposed in the previous section is determined by the two steps needed to compute the estimate:
\begin{enumerate}
	\item Drawing $M$ frequency point samples from the PDF $p$.
	\item Evaluating the transfer function $f(\mu; s)$ at the $M$ frequency point samples.
\end{enumerate}
In determining computational complexity, note that each frequency point sample is drawn independently. Thus, there is a high potential for parallel execution. The same holds for evaluating the transfer function $f(\mu; s)$. 
\subsubsection{Drawing frequency point samples}
Sampling from a PDF can be achieved using various techniques, depending on the type of PDF. In most computational environments, we have access to a routine that draws a uniform random number between 0 and 1. We will now show how this routine can be used to draw random numbers from more general distributions.

As a first example, consider a \emph{log-uniform} distribution between frequencies $\omega_{\text{min}} \in (0, \infty)$ and $\omega_{\text{max}} \in [\omega_{\text{min}}, \infty)$. By drawing a uniform random number $\alpha$ between 0 and 1 and scaling it to lie in the range $[\log_{10}{\omega_{\text{min}}}, \log_{10}{\omega_{\text{max}}}]$, the log-uniformly distributed number is given as $10^{\alpha}$.

For more general distributions, the inversion method \cite{Devroye2006} is appropriate if the inverse of the cumulative density function (CDF) - denoted by $F^{-1}(\omega)$ - can be efficiently evaluated. In that case, we can draw a uniform random number $\alpha$ between 0 and 1 and then solve the equation $F^{-1}(\omega) = \alpha$ for $\omega$.

In the sequel, we assume that frequency point samples can be drawn by applying elementary operations to uniformly drawn random numbers (such as log-uniform distributions) or that evaluation of the inverse CDF is relatively cheap and does not scale with $n$. As a result, this step in the computation takes $\mathcal{O}(M)$ time on a non-parallel implementation, with a further reduction if parallel computing resources are available.

\subsubsection{Evaluating the frequency response of $f$}
The system $f(\mu)$ (\ref{eqn:f}) is composed of the large-scale dynamical system that should be optimized ($G(\mu)$) and its parameter gradient systems ($\frac{\partial G}{\partial \mu_j}$). Thus, if we can efficiently evaluate frequency samples of these systems, we can efficiently evaluate frequency samples of $f(\mu)$. In scenarios where the parameter-dependent system originates from the semi-discretization of a Partial Differential Equation (PDE), we have access to a sparse state-space description of these systems (for example, by extracting it from the finite element analysis software) in the form
\begin{equation}
	\begin{aligned}
		E \dot{x} & = A x + B u, \\
				y & = C x.
	\end{aligned}
\end{equation}
The frequency response of such a state-space representation at frequency $i\omega_m$ is given by $C (i\omega_m E - A)^{-1} B$, which can be efficiently calculated in $\mathcal{O}(n)$ time by using linear solvers designed for sparse systems \cite{Saad2003}. As a result, this step in the computation takes $\mathcal{O}(M n)$ time on a non-parallel implementation.

\subsubsection{Total computational complexity}
Summing the computational complexity of the two stages, it can be seen that the computational complexity of the gradient estimator is $\mathcal{O}(M ) + \mathcal{O}(M n) = \mathcal{O}(M n)$. Note that the proposed estimator is amenable to parallelization (both in drawing the $M$ samples from the PDF $p$ and calculating the frequency response of $f(\mu)$ at those samples).

\section{Numerical examples} \label{sct:numerical_examples}
In this section, the proposed method is evaluated on two numerical examples. In Section \ref{sct:observer_design} we consider a fixed-order observer design problem for a large-scale thermal model and in Section \ref{sct:pid_tuning} we consider a controller tuning problem for an infinite-dimensional mechanics model. We used Algorithm \ref{alg:sh2opt} implemented in MATLAB version R2023b.

\subsection{Fixed-order observer design} \label{sct:observer_design}
\begin{figure}
	\centering
	\includegraphics[width=0.35\textwidth]{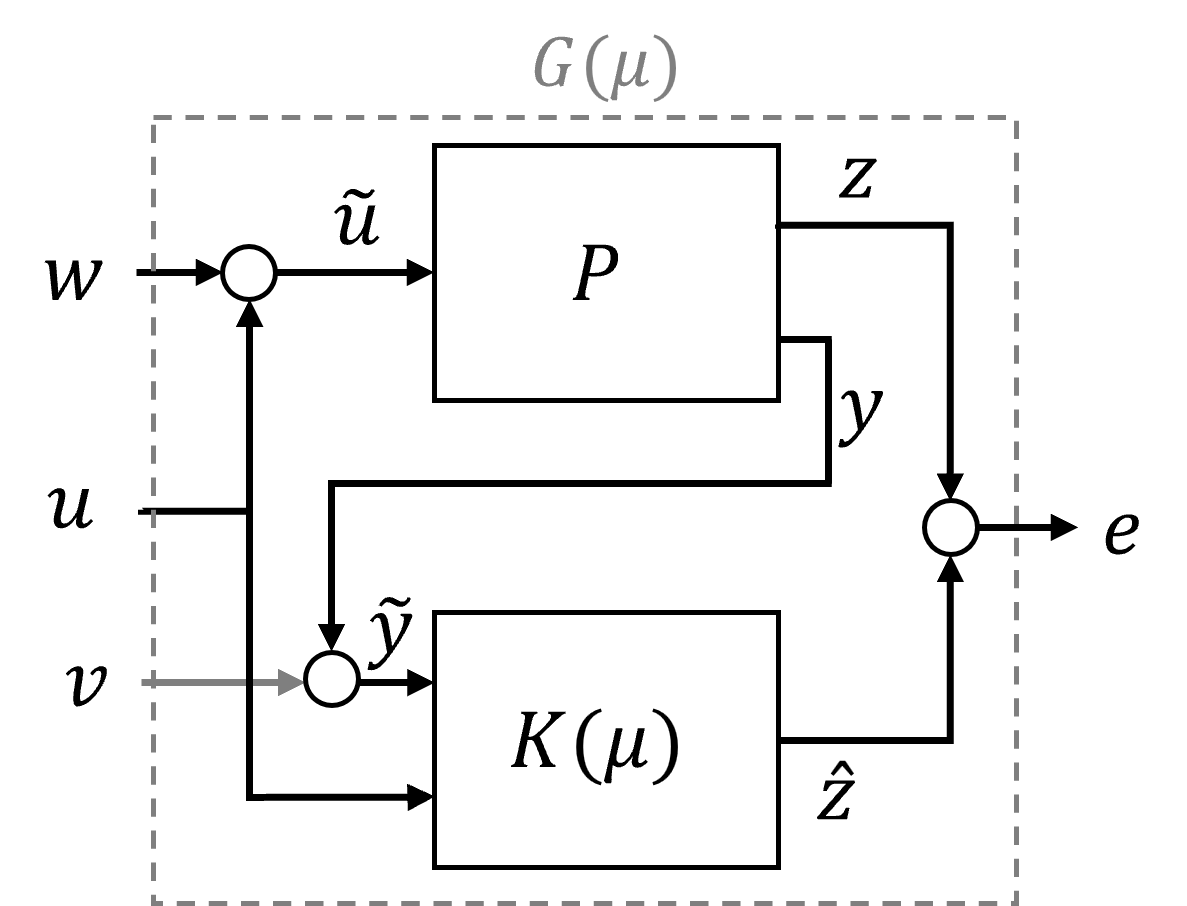}
	\caption{Fixed-order observer design problem considered in Section \ref{sct:observer_design}.}
	\label{fig:observer_block_diagram}
\end{figure}

In this section, we consider a fixed-order observer design problem. The plant model is a semi-discretized thermal model of a chip cooling system \cite{morwiki_convection}:
\begin{equation}
	P : \begin{cases}
		E \dot{x} & = A x + B (u + w), \\
			   z & = C_z x, \\
			   y & = C_y x, \\
			   \tilde{y} &= y + v.
	\end{cases}
	\label{eqn:observer_design_P}
\end{equation}
with state vector $x(t) \in \mathbb{R}^{20,082}$, deterministic input $u(t) \in \mathbb{R}$, measured output $\tilde{y}(t) \in \mathbb{R}$ and the output-to-be-observed $z(t) \in \mathbb{R}$, which is not measured. The white noise signals $w(t) \in \mathbb{R}$ and $v(t) \in \mathbb{R}$ represent the uncertainty. We are interested in developing an observer producing an estimate $\hat{z}$ of $z$ that minimizes the $\mathcal{H}_2$ norm of the transfer function $G(\mu)$ mapping $(w, u, v)$ to the estimation error $e := z - \hat{z}$ (see Figure \ref{fig:observer_block_diagram}). It is well-known that the optimal solution corresponds to the Kalman filter \cite{Kalman1960}. However, the Kalman filter is of the same order of the plant and computing it incurs a computational complexity that scales cubically with the number of states of the plant. Motivated by this, we seek instead a locally optimal $\mathcal{H}_2$ \emph{fixed-order} observer of the form
\begin{equation}
	K(\mu): \begin{cases}
		\dot{q} &= A_q(\mu) q + B_q(\mu) \begin{bmatrix} u \\ \tilde{y} \end{bmatrix}, \\
		 	 \hat{z} &= C_q(\mu) q,
			\end{cases}
\end{equation}
with $q(t) \in \mathbb{R}^{2}$ and $A_q(\mu)$, $B_q(\mu)$, $C_q(\mu)$ appropriately-sized matrix-valued functions in $\mu$ such that
\begin{equation}
	\mathrm{vec}\begin{bmatrix} A_q(\mu) & B_q(\mu) & C_q(\mu)^{\mathrm{T}} \end{bmatrix} = \mu \in \mathbb{R}^{10}.
\end{equation}
To initialize Algorithm \ref{alg:sh2opt}, we reduce the plant (\ref{eqn:observer_design_P}) to order $r = 2$ by IRKA \cite{Gugercin2008} and synthesize a Kalman filter for the reduced-order plant. From the state-space coefficients of this Kalman filter, we extract $\mu_0$. We also scale $C_z$ in (\ref{eqn:observer_design_P}) to normalize the estimation error with respect to the absence of an observer ($\hat{z} = 0$), i.e., $\hnorm{G(0)} = 1$. We approximate the $\mathcal{H}_2$ norm required for this normalization step and the subsequent evaluation of the results using the sparse Lyapunov equation solver from version 3 of the M-M.E.S.S. toolbox \cite{SaaKB21-mmess-3.0}. The step size policy is as follows:
\begin{equation}
	\alpha_k = \begin{cases}
			10^{-4} & \text{for} \quad k \leq 20, \\
			\frac{1}{2} 10^{-4} & \text{for} \quad 20 < k \leq 40, \\
			\frac{1}{4} 10^{-4} & \text{for} \quad 40 < k \leq 60, \\
			\frac{61}{k} \frac{1}{8} 10^{-4} & \text{for} \quad 60 < k.
		\end{cases}
\end{equation}
Note that this step size policy is $\alpha_k = \Theta(k^{-1})$ which is $\ell_2$ summable. Hence, this step size is in line with the findings from the stability and convergence analysis from Section \ref{sct:convergence}.

\begin{figure}
	\centering
	\includegraphics[width=0.49\textwidth]{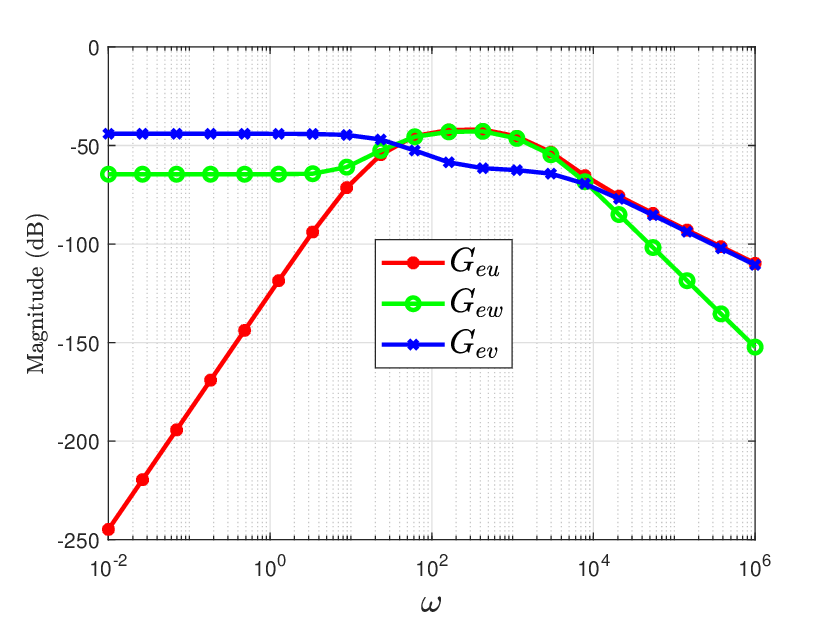}
	\caption{Bode magnitude diagram of $G(\mu_0)$.}
	\label{fig:observer_design_bode}
\end{figure}

\begin{figure}
	\centering
	\includegraphics[width=0.49\textwidth]{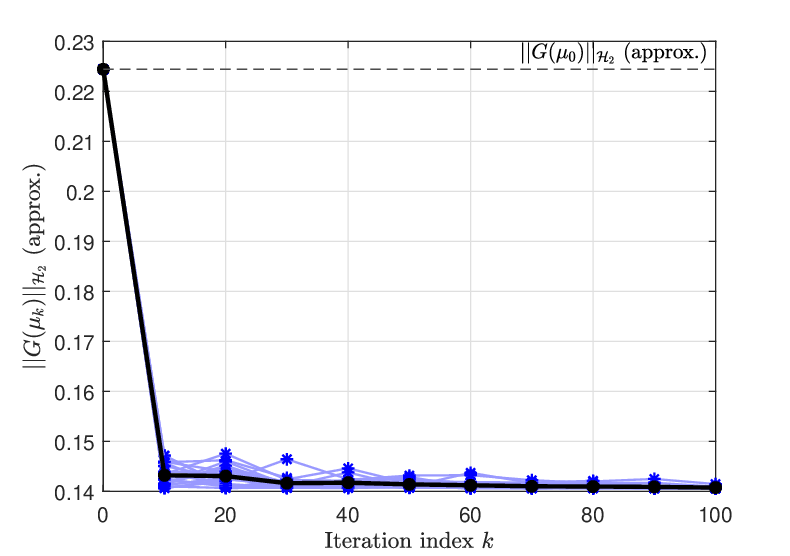}
	\caption{Approximate $\mathcal{H}_2$ norm of $G(\mu_k)$. Thin lines represent separate runs and the thick line is the mean.}
	\label{fig:observer_design_results}
\end{figure}

Next, an appropriate sampling PDF was determined. This was done by first estimating the relevant frequency range by plotting the Bode magnitude diagram of $G(\mu_0)$ at 20 frequency points on a logarithmic scale as shown in Figure \ref{fig:observer_design_bode}. A loguniform sampling PDF with a frequency range from $10^{-2}$ to $10^{6}$ rad/s is selected. We set $M = 10$, i.e., we take 10 frequency samples distributed according to the PDF in each iteration. Since the method is stochastic, we run 20 trials using the same initialisation, step size and sampling PDF. The approximate $\mathcal{H}_2$ norm of all runs for every 10 iterations is given in Figure \ref{fig:observer_design_results}. In every run, the estimation error of the fixed-order observed is improved from the initial Kalman filter based on the reduced-order model ($\hnorm{G(\mu_0)} = 0.22$) to the optimized observer $\hnorm{G(\mu_{100})} = 0.14$). It can also be observed that the variance between runs becomes smaller as the step size is decreased and all runs appear to converge towards a local minimizer.

\begin{figure}
	\centering
	\includegraphics[width=0.49\textwidth]{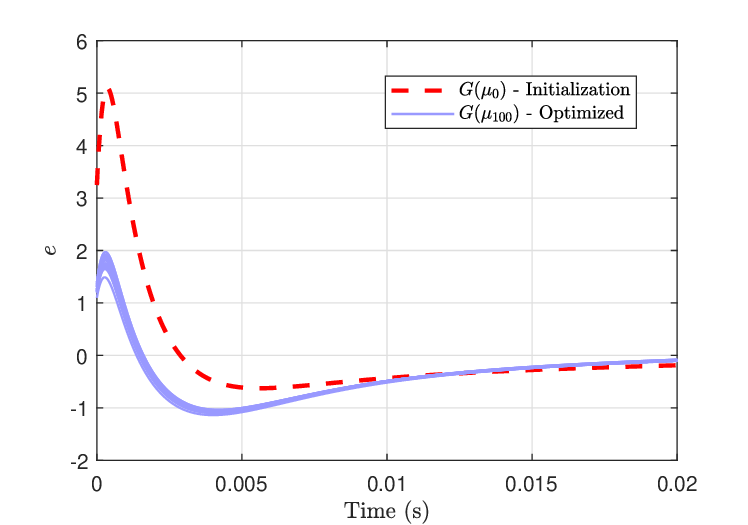}
	\caption{Comparison of impulse response of $G_{eu}(\mu)$ between initial and all optimized observers.}
	\label{fig:observer_design_simulation}
\end{figure}

To illustrate the improved observer performance, the impulse response of the channel $G_{eu}(\mu)$ is shown in Figure \ref{fig:observer_design_simulation} for the initial ($\mu = \mu_0$) and optimized ($\mu = \mu_{100}$) observers. All optimized observers show a similar dynamic response. The peak in the transient estimation error seen for $\mu = \mu_0$ is attenuated significantly by the optimization in all runs from 5 to 2 and lower.

\subsection{PD controller tuning for an infinite-dimensional system} \label{sct:pid_tuning}
\begin{figure}
	\centering
	\includegraphics[width=0.35\textwidth]{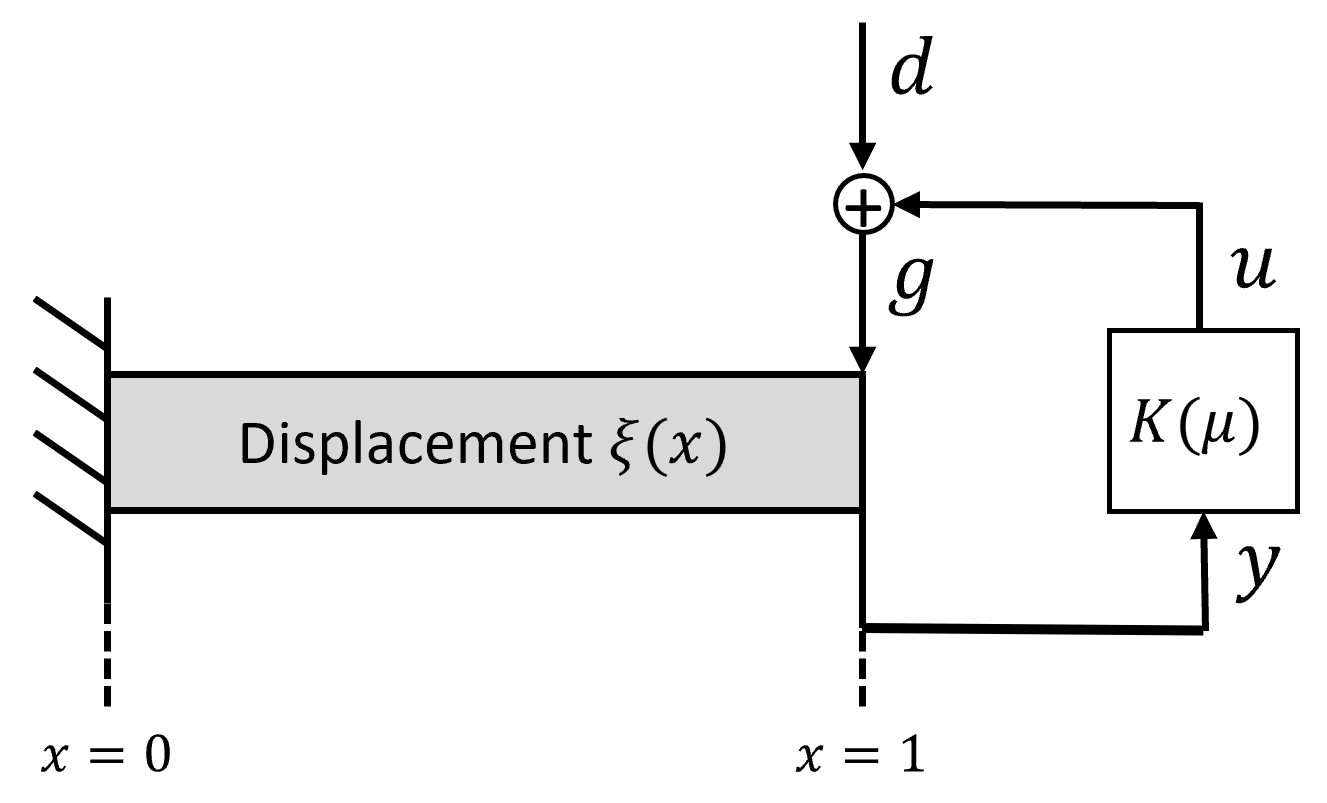}
	\caption{Block diagram of the PD tuning problem considered in Section \ref{sct:pid_tuning}. $K(\mu)$ represents the parametrized PD controller. The control cost is captured by $G(\mu): d \rightarrow (y, u)$.}
	\label{fig:pid_tuning_block_diagram}
\end{figure}

We consider a cantilever beam on a domain $x \in [0, 1]$ for which the displacement $\xi(x)$ is governed by the 1-dimensional damped wave equation:
\begin{equation}
	\frac{\partial^2 \xi}{\partial t^2} + c \frac{\partial \xi}{\partial t} - \frac{\partial^2 \xi}{\partial x^2} = 0,
	\label{eqn:pid_tuning_pde}
\end{equation}
with $c = 0.25$. At the fixed end ($x = 0$) we have $\xi(0) = 0$. At the free end ($x = 1$) we measure the displacement $y := \xi(1)$ and denote by $g$ the net force. A controller $K(\mu)$ measures $y$ and outputs a control force $u$ as displayed in Figure \ref{fig:pid_tuning_block_diagram}. The control cost to be minimised encoded in $G(\mu)$ consists of two input-to-output channels: 1) the transfer from $d$ to $y$; and 2) the transfer from $d$ to $u$. Note that our method enables the user to easily add weighting filters to steer performance to application-specific relevant frequency bands. For the controller $K(\mu)$, we use a PD controller with filtered derivative term parametrized as follows:
\begin{equation}
	K(\mu; s) = \left[ \mu \right]_1 + \left[ \mu \right]_2 \frac{T_F s}{T_F s + 1},
\end{equation}
with $T_F = 10^{-2}$.

In this setting, an advantage of the proposed SH2OPT method is that we do not need a (finite-dimensional) state-space realization of $G(\mu)$. Thus, semi-discretization of (\ref{eqn:pid_tuning_pde}) is not required (which would introduce an approximation error). Instead, we can directly work with the analytical solution of the wave equation (\ref{eqn:pid_tuning_pde}) in the frequency domain based on the general solution of the form
\begin{equation}
	\Xi(x) = c_1 e^{\phi x} + c_2 e^{-\phi x},
\end{equation}
where $\Xi(x) \in \mathbb{C}$ denotes the phasor associated with $\xi$ (i.e., $\xi(x, t) = \mathrm{Re}(\Xi(x) e^{i\omega t})$), $c_1, c_2 \in \mathbb{C}$ are coefficients that depend on the boundary conditions and $\phi = \sqrt{ic - \omega^2}$.
Application of the boundary condition at $x = 0$ gives $c_2 = -c_1$. At $x = 1$ we have the prescribed force $g$ leading to
\begin{equation}
	c_1 = \frac{g}{\phi(e^{\phi} + e^{\phi})}.
\end{equation}
Based on the above, we obtain the following expression for $G(\mu)$:
\begin{equation}
	G(\mu; i\omega) = \begin{bmatrix} 1 \\ K(\mu; i\omega) \end{bmatrix} \frac{\Phi(\omega)}{1 - \Phi(\omega) K(\mu; i\omega)},
	\label{eqn:pid_tuning_G}
\end{equation}
with
\begin{equation}
	\Phi(\omega) := \frac{e^{\phi(\omega)} - e^{-\phi(\omega)}}{\phi(\omega) \left( e^{\phi(\omega)} + e^{-\phi(\omega)} \right)}.
\end{equation}
The analytical expression (\ref{eqn:pid_tuning_G}) allows us to efficiently draw the frequency samples needed in Algorithm \ref{alg:sh2opt}.

\begin{figure}
	\centering
	\includegraphics[width=0.49\textwidth]{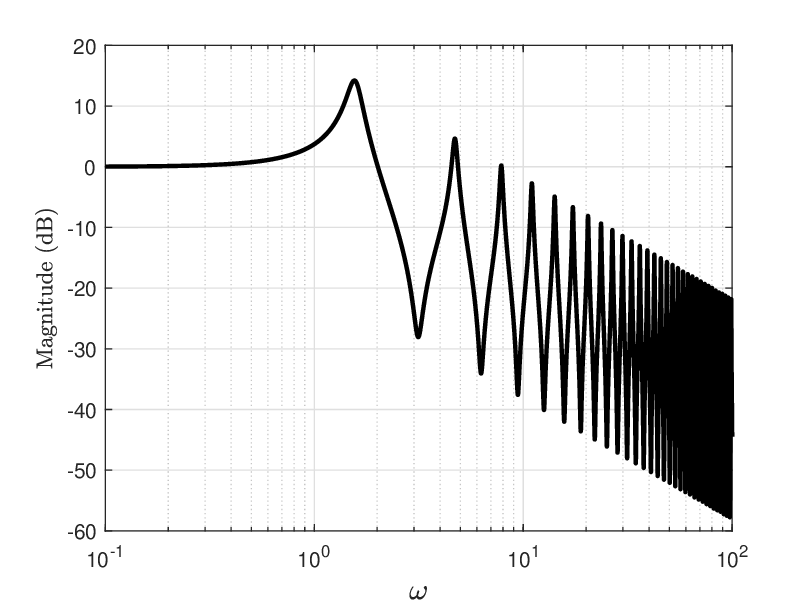}
	\caption{Bode magnitude diagram of $G(\mu_0)$.}
	\label{fig:pid_tuning_bodemag}
\end{figure}

\begin{figure}
	\centering
	\includegraphics[width=0.49\textwidth]{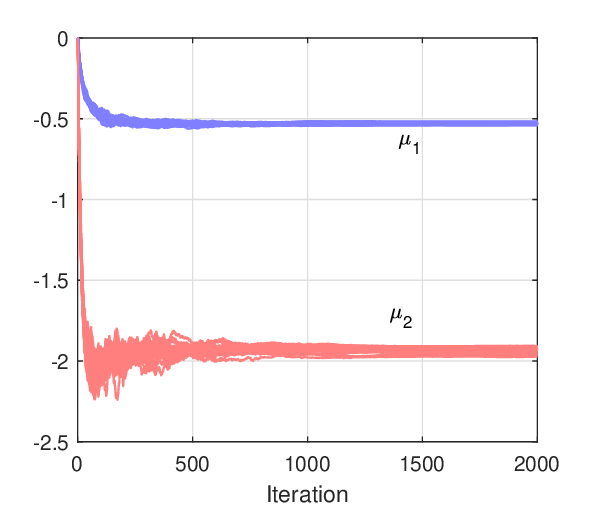}
	\caption{Evolution of the PD parameters. Each line represents a separate run.}
	\label{fig:pid_tuning_results_parameters}
\end{figure}

We initialize Algorithm \ref{alg:sh2opt} using the open-loop scenario ($\mu_0 = 0$). Note that $G(\mu_0)$ is guaranteed to be stable due to the positive damping factor $c = 0.25$. Furthermore, $\muSpace$ is selected such that the PD parameters are non-positive to obtain closed-loop stability.
The Bode magnitude diagram of $G(\mu_0)$ is shown in Figure \ref{fig:pid_tuning_bodemag}. Based on this diagram, the PDF for this numerical experiment is chosen as a loguniform distribution between $10^{-2}$ and $10^4$ rad/s. At each iteration, $M = 1000$ frequency samples are drawn from this PDF to estimate the gradient. The step size is initialized at $\alpha_0 = 10^{-2}$ and halved every 200 iterations. The experiment is repeated 20 times and the evolution of the parameters is plotted for each trial in Figure \ref{fig:pid_tuning_results_parameters}. A rapid decrease in both parameters values is observed for all runs during the first 50 iterations. Afterwards, the spread of parameter values between trails decreases during the first 1000 iterations, after which no significant decrease in spread is observed. The remaining spread we observe in the final parameter values across the runs could be caused by a low sensitivity of the $\mathcal{H}_2$ norm near those values. To investigate this further, an approximation of the achieved $\mathcal{H}_2$ norm was computed by a semi-discretised approximation $\tilde{G}(\mu) \approx G(\mu)$ using the finite-difference method with $n = 400$ states for the final parameter value of all runs. The computed $\mathcal{H}_2$ norm was scaled by the $\mathcal{H}_2$ norm for the initial parameter, i.e., $\hnorm{\tilde{G}(\mu_0)} = 1$. The results for the 20 runs were as follows. The average $\mathcal{H}_2$ norm was $8.29 \cdot 10^{-1}$ and the standard deviation was $3.89 \cdot 10^{-5}$. The low standard deviation indeed supports the statement that a relatively large spread of the final parameter values is caused by a low sensitivity of the $\mathcal{H}_2$ norm near those values.

\begin{figure}
	\centering
	\includegraphics[width=0.49\textwidth]{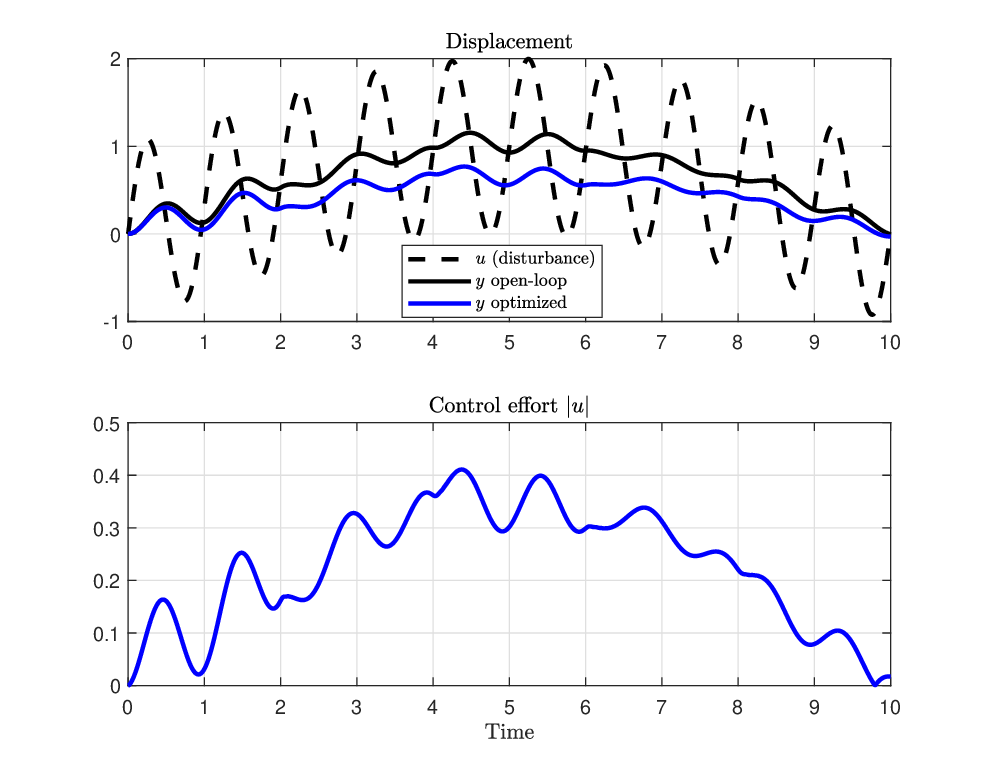}
	\caption{Time response of the approximate closed-loop system $\tilde{G}(\mu)$ excited by a disturbance.}
	\label{fig:pid_tuning_simulation}
\end{figure}

The transient performance of the optimized PD controller obtained in the last run is studied by a simulation using the approximate system model $\tilde{G}(\mu)$. The system is simulated for $t \in [0, 10]$. The input disturbance is $u(t) = \sin(2\pi t) + \sin(2\pi 10^{-1} t)$. The response in terms of the actuation effort $u$ and displacement $y$ is plotted in Figure \ref{fig:pid_tuning_simulation}. The open-loop scenario ($\mu = \mu_0$) and the optimized controller are compared. It can be seen that the controller attenuates the displacement $y$. Since the control cost encoded by $G(\mu)$ includes also the control effort $u$, the resulting performance balances the attenuation in displacement and required control effort. In practical scenarios, the user can influence this balance by selecting (frequency-dependent) weighting filters.

\section{Concluding remarks} \label{sct:concluding_remarks}
In this paper, we have proposed an SGD optimization scheme (SH2OPT) for large-scale dynamical systems that is computationally feasible and, in contrast to existing techniques, provides (probabilistic) stability and convergence guarantees for a wide class of differentiable parametrized dynamical systems. In this way, we have opened the way to transfer the theoretical and practical success of SGD as observed in large-scale static optimization problems to dynamical optimization problems.

SH2OPT was tested on numerical examples in the areas of fixed-order observer design and controller synthesis for infinite-dimensional plants. These numerical examples demonstrated the potential of the approach.

Future improvements to the proposed method can be achieved by variance reduction schemes that adaptively update the sampling distribution. This can be achieved, for example, by interpolating frequency domain samples of $f(\mu; i\omega)$ as discussed in Section \ref{sct:pj}. However, care must be taken to ensure guarantees on stability preservation remain applicable. Besides variance reduction, additional improvements can be expected by considering the improvements that have been made to SGD in static problems, such as formulating a second-order scheme similar to \cite{Bordes2009} or considering implicit SGD updates as in  \cite{Toulis2017} to decrease sensitivity to settings such as the step size.

\begin{ack}
This work has been funded in part by ITEA under the COMPAS project (ITEA project 19037).
\end{ack}

\bibliographystyle{plain}
\bibliography{library}

\begin{appendix}
\section{Stability guarantees} \label{sct:boundedness_of_iterates}
In this section, we prove Theorem \ref{thm:stability}. In the proof, we will make use of the following lemma.
\begin{lem} \label{lem:norm_xi}
Let $x, y \in \mathbb{R}^{n}$. Then,
\begin{equation}
 \max_{\xi \in [0, 1]} \norm{x} \norm{\xi x + y} \leq 2 \left( \norm{x}^2 + \norm{y}^2 \right)
 \label{eqn:norm_xi}
\end{equation}
\end{lem}
\begin{pf}
Consider first the case that $x^{\mathrm{T}} y \geq 0$. Then,
\begin{equation}
\begin{aligned}
\max_{\xi \in [0, 1]} \norm{x} \norm{\xi x + y} \leq \norm{x}^2 + x^\mathrm{T} y.
\end{aligned}
\end{equation}
Thus, (\ref{eqn:norm_xi}) holds if we can show that
\begin{equation}
	\norm{x}^2 + x^{\mathrm{T}} y \leq 2\left( \norm{x}^2 + \norm{y}^2 \right).
\end{equation}
The following chain of inequalities shows that this is indeed the case.
\begin{equation}
\begin{aligned}
\norm{x}^2 + x^{\mathrm{T}} y & \leq 2\left( \norm{x}^2 + \norm{y}^2 \right) \\
\Rightarrow x^\mathrm{T} y & \leq \norm{x}^2 + 2 \norm{y}^2 \\
\Rightarrow x^\mathrm{T} y & \leq \norm{y} + \norm{x - y}^2 + 2 x^\mathrm{T} y \\
\Rightarrow 0 & \leq \norm{y} + \norm{x - y}^2 + x^\mathrm{T} y.
\end{aligned}
\end{equation}
Next, consider the case that $x^\mathrm{T} y < 0$. We have:
\begin{equation}
\max_{\xi \in [0, 1]} \norm{x} \norm{\xi x + y} = \norm{x} \norm{y} \leq -x^\mathrm{T} y.
\end{equation}
Thus, (\ref{eqn:norm_xi}) holds if we can show that
\begin{equation}
	-x^\mathrm{T} y \leq 2\left( \norm{x}^2 + \norm{y}^2 \right).
\end{equation}
The following chain of inequalities shows that, again, this is the case.
\begin{equation}
\begin{aligned}
-x^\mathrm{T} y & \leq 2 \left( \norm{x}^2 + \norm{y}^2 \right) \\
\Rightarrow 0 & \leq \frac{3}{2} \norm{x}^2 + \frac{3}{2} \norm{y}^2 + \frac{1}{2} \norm{x}^2 + x^\mathrm{T} y + \frac{1}{2} \norm{y}^2 \\
\Rightarrow 0 & \leq \frac{3}{2} (\norm{x}^2 + \norm{y}^2) + \frac{1}{2} \norm{x + y}^2
\end{aligned}
\end{equation}
 $\blacksquare$
\end{pf}
The proof of Theorem \ref{thm:stability} is divided in 3 steps:
\begin{enumerate}
	\item The change in the cost between subsequent iterations, i.e., $c(\mu_{k+1}) - c(\mu_{k})$, is decomposed in two additive terms: a deterministic term $D_k$ and a stochastic term $S_k$. Using Lipschitz smoothness of $c$, a condition on the step size $\alpha_k$ is derived which guarantees non-positivity of $D_k$.
	\item The stochastic term $S_k$ is further decomposed in two additive terms which are linear, respectively quadratic in the step size $\alpha_k$.
	\item Following the line of reasoning in the proof of \cite[Theorem 4]{Mertikopoulos2020}, the decomposition of $S_k$ derived in the preceding step guarantees boundedness of the series $\sum_{k=0}^{\infty} S_k$ with probability at least $1 - \delta$ provided the step size $\alpha_k$ decreases sufficiently rapidly.
\end{enumerate}

\textbf{Step 1} $\cdot$ Let $\mu_{k} \in \muSpaceStab$. If also the next iterand $\mu_{k+1} \in \muSpaceStab$ (we will discuss at the end of this step how to ensure this is a reasonable assumption), then using (\ref{eqn:sgd_static}) the change in cost between these iterates can be decomposed as
\begin{equation}
\begin{aligned}
 c(\mu_{k+1}) - c(\mu_k) & = \int_{\mu_k}^{\mu_{k+1}} \nabla c(\tau)^\mathrm{T} d\tau \\
	& = \underbrace{\int_{\mu_k}^{\mu_k - \alpha_k \nabla c(\mu_k)} \nabla c(\tau)^\mathrm{T} d\tau}_{D_k} \\
	& + \underbrace{\int_{\mu_k - \alpha_k \nabla c(\mu_k)}^{\mu_k -\alpha_k (\nabla c(\mu_k) + Z_k)} \nabla c(\tau)^\mathrm{T} d\tau}_{S_k}.
\end{aligned}
	\label{eqn:Dk_Sk_decomposition}
\end{equation}
The decomposition splits $c(\mu_{k+1}) - c(\mu_k)$ into a deterministic component $D_k$ and a stochastic component $S_k$. Using the mean value theorem, the deterministic term $D_k$ can be written as
\begin{equation}
D_k = -\alpha_k \nabla c(\mu_k)^{\mathrm{T}} \nabla c(\mu_k - \xi \alpha_k \nabla c(\mu_k))
\label{eqn:Dk}
\end{equation}
for some $\xi \in [0, 1]$. Since $\alpha_k \geq 0$, $D_k$ in (\ref{eqn:Dk}) is guaranteed to be non-positive if
\begin{equation}
	\norm{\nabla c(\mu_k) - \nabla c(\mu_k - \xi \alpha_k \nabla c(\mu_k))} \leq \norm{\nabla c(\mu_k)},
	\label{eqn:lipschitz_bound}
\end{equation}
which holds trivially if $\nabla c(\mu_k) = 0$. Otherwise, applying Lipschitz smoothness of $c(\mu_k)$ gives the following sufficient condition for (\ref{eqn:lipschitz_bound}):
\begin{equation}
	L \norm{ \xi \alpha_k \nabla c(\mu_k) } \leq \norm{\nabla c(\mu_k}.
	\label{eqn:Dk_nonpositive_condition}
\end{equation}
Noting that $\xi \in [0, 1]$, the inequality (\ref{eqn:Dk_nonpositive_condition}) is satisfied if $\alpha_k$ satisfies
\begin{equation}
	\alpha_k \leq L^{-1},
	\label{eqn:alpha_k_bound}
\end{equation}
which is guaranteed under the conditions of the theorem.

To guarantee stable iterands for all $k = 0, 1, ...$, we must show that the cumulative effect of the deterministic ($D_k$) and stochastic ($S_k$) terms are bounded from above. The bound on $\alpha_k$ given by (\ref{eqn:alpha_k_bound}) guarantees that $D_k \leq 0$, hence also guaranteeing that the cumulative effect of $D_k$ is bounded from above by 0. In the remaining steps, we are therefore left with determining that, indeed, the cumulative effect of $S_k$ is also bounded from above with high probability.

\textbf{Step 2} $\cdot$ If (\ref{eqn:alpha_k_bound}) is satisfied (and hence $D_k \leq 0$), the cost of the k-th iteration is bounded by the sum of the stochastic components of (\ref{eqn:Dk_Sk_decomposition}):
\begin{equation}
c(\mu_k) \leq c(\mu_0) + \sum_{i=0}^{k-1} S_i.
\label{eqn:ck_bound}
\end{equation}
Similar to (\ref{eqn:Dk}), we can use the mean value theorem to write the random variable $S_k$ as
\begin{equation}
S_k = -\alpha_k Z_k^{\mathrm{T}} \nabla c(\mu_k - \alpha_k \nabla c(\mu_k) - \xi \alpha_k Z_k)
\end{equation}
for some $\xi \in [0, 1]$. Subtracting and adding $\alpha_k Z_k^{\mathrm{T}} \nabla c(\mu_k)$ and applying Lipschitz smoothness of $c(\mu_k)$ leads to the following bound:
\begin{equation}
\begin{aligned}
S_k & \leq -\alpha_k Z_k^{\mathrm{T}} \nabla c(\mu_k) + L \alpha_k^2 \norm{Z_k} \norm{\nabla c(\mu_k) + \xi Z_k} \\
& \leq -\alpha_k Z_k^{\mathrm{T}} \nabla c(\mu_k) + 2 L \alpha_k^2 \left( \norm{\nabla c(\mu_k)}^2 + \norm{Z_k}^2 \right),
\end{aligned}
\label{eqn:Sk_bound}
\end{equation}
where the second inequality follows from Lemma \ref{lem:norm_xi}. Inequality (\ref{eqn:Sk_bound}) gives us a bound on $S_k$ as the sum of two terms. The first term is zero-mean (since $\mathbb{E}[Z_k] = 0$) and thus has (in expectation) a net zero contribution to the iterates. Although the second term is non-negative, we will show in step 3 of the proof that, since it is quadratic in $\alpha_k$, it decays sufficiently fast if condition (\ref{eqn:alpha_sum_bound}) in the theorem is satisfied. To simplify notation in the sequel, introduce $\phi_i := -Z_i^\mathrm{T} \nabla c(\mu_i)$ and $\psi_i := 2L \left( \norm{\nabla c(\mu_i)}^2 + \norm{Z_i}^2 \right)$. Then, we can write (\ref{eqn:Sk_bound}) compactly as:
\begin{equation}
S_k \leq \alpha_k \phi_k + \alpha_k^2 \psi_k.
\label{eqn:Sk_bound_compact}
\end{equation}
The bound (\ref{eqn:ck_bound}) can then be written as
\begin{equation}
c(\mu_k) \leq c(\mu_0) + M_{k-1} + N_{k-1},
\label{eqn:ck_bound_telescoped}
\end{equation}
with
\begin{equation}
M_{k-1} := \sum_{i=0}^{k-1} \alpha_i \phi_i
\end{equation}and
\begin{equation}
N_{k-1} := \sum_{i=0}^{k-1} \alpha_i^2 \psi_i.
\label{eqn:Nk}
\end{equation}

\textbf{Step 3} $\cdot$ We follow the steps of the proof of \cite[Theorem 4]{Mertikopoulos2020} to arrive at the final result.
We first introduce the event $E_k$ that all iterates are contained in the sublevel set $\mathcal{U} = \{ \mu | c(\mu) \leq c(\mu_0) + \sqrt{\varepsilon} + \varepsilon \}$:
\begin{equation}
	E_k := \{ \mu_i \in \mathcal{U} \text{ for all } i = 0, 1, ..., k-1 \}.
\end{equation}
We also define the error quantity
\begin{equation}
R_k := M_k^2 + N_k
\label{eqn:Rk}
\end{equation}
and the associated event $H_k$ that $R_k$ is not larger than $\varepsilon$ for the first $k$ iterations, i.e.:
\begin{equation}
	H_k := \{ R_i \leq \varepsilon \text{ for all } i = 0, 1, ..., k-1 \}.
\end{equation}
The following lemma connects the preservation of stability (event $E_k$) with a bound on the stochastic error $R_k$ (event $H_k$).
\begin{lem} \label{lem:Ek_Hk}
For all $k = 1, 2, 3, ...$ we have $H_k \subseteq E_{k+1}$.
\end{lem}
\begin{pf}
Fix a realization of $H_k$ so that $R_i \leq \varepsilon$ for all $i = 0, 1, ..., k-1$. Then, we need to show that $E_{k+1}$ occurs, which we do by combining (\ref{eqn:ck_bound_telescoped}) and (\ref{eqn:Rk}) to arrive at the following bound:
\begin{equation}
\begin{aligned}
c(\mu_i) & \leq c(\mu_0) + M_{i-1} + N_{i-1} \\
	& \leq c(\mu_0) + \sqrt{R_{i-1}} + R_{i-1} \\
	& \leq c(\mu_0) + \sqrt{ \varepsilon } + \varepsilon,
\end{aligned}
\end{equation}
where the second inequality follows from the definition (\ref{eqn:Rk}) and the fact that $N_k \geq 0$ by (\ref{eqn:Nk}).
The above inequality holds for all $i = 1, 2, ..., k$. Hence, $E_{k+1}$ occurs. $\blacksquare$
\end{pf}
In the last part of the proof, we will show that the probability that $H_k$ occurs for all iterations $k = 1, 2, 3, ...$ can be moved arbitrarily close to 1. To this end, we introduce another event $\widetilde{H}_k$ as follows:
\begin{equation}
\begin{aligned}
\widetilde{H}_k & := H_{k-1} \setminus H_{k} \\
& = \{ R_i \leq \varepsilon \text{ for all } i = 0, 1, ..., k - 2 \text{ and } R_{k-1} > \varepsilon \}.
\end{aligned}
\end{equation}
The event $\widetilde{H}_k$ encodes the occurrence of a large contribution of the stochastic term at the k-th iteration. Furthermore, introduce $\widetilde{R}_k = R_k \indicator_{H_k}$ with $\indicator_{H_k}$ the indicator function, i.e.,
\begin{equation}
	\indicator_{H_k} = \begin{cases} 
			0 \text{ if } H_k \text{ does not occur}, \\
			1 \text{ if } H_k \text{ occurs}.
		\end{cases}
\end{equation}
The following lemma bounds the growth of $\widetilde{R}_k$ (in expectation).
\begin{lem} \label{lem:Hk_tilde}
\begin{equation}
	\mathbb{E}[\widetilde{R}_{k+1}] \leq \mathbb{E}[\widetilde{R}_{k}] + R_{*} \alpha_k^2 - \varepsilon \mathbb{P}[\widetilde{H}_{k+1}],
	\label{eqn:Rk_lemma}
\end{equation}
with $R_{*} = K^2 \sigma^2 + 2L(K^2 + \sigma^2)$.
\end{lem}
\begin{pf}
We start by decomposing $\widetilde{R}_k$ as
\begin{equation}
\begin{aligned}
	\widetilde{R}_{k} & = R_{k} \indicator_{H_{k}} \\
		& = R_{k-1} \indicator_{H_{k}} + (R_{k} - R_{k-1})\indicator_{H_{k}} \\
		& = R_{k-1} \indicator_{H_{k-1}} - R_{k-1} \indicator_{\widetilde{H}_{k}} + (R_{k} - R_{k-1})\indicator_{H_{k}} \\
		& = \widetilde{R}_{k-1} + (R_{k} - R_{k-1})\indicator_{H_{k}} - R_{k-1} \indicator_{\widetilde{H}_{k}}.
\end{aligned}
\label{eqn:Rk_decomposed}
\end{equation}
with the third equality following from the fact that $H_k = H_{k-1} \setminus \widetilde{H}_k$ and hence $\indicator_{H_k} = \indicator_{H_{k-1}} - \indicator_{\widetilde{H}_k}$. We will show that taking the expectation of the terms in (\ref{eqn:Rk_decomposed}) gives the result in (\ref{eqn:Rk_lemma}). We begin with the term $(R_k - R_{k-1})\indicator_{H_k}$ by writing the following recursive relation for $R_k$:
\begin{equation}
\begin{aligned}
R_{k+1} & = (M_{k} + \alpha_{k} \phi_{k})^2 + N_{k} + \alpha_{k}^2 \psi_{k} \\
	& = M_{k}^2 + N_{k} + 2 M_k \alpha_k \phi_k + \alpha_k^2 \phi_k^2 + \alpha_k^2 \psi_k \\
	& = R_k + 2M_k \alpha_k \phi_k + \alpha_k^2 \phi_k^2 + \alpha_k^2 \psi_k.
\end{aligned}
\label{eqn:Rk_recursive}
\end{equation}
Next, we use (\ref{eqn:Rk_recursive}) to write
\begin{subequations}
\begin{align}
	\mathbb{E}[\left( R_{k+1} - R_{k} \right) \indicator_{H_k}] & = 2 \alpha_k \mathbb{E}[M_{k} \phi_k \indicator_{H_k}] \label{eqn:Rk_term1} \\
		& + \alpha_k^2 \mathbb{E}[\phi_k^2 \indicator_{H_k}] \label{eqn:Rk_term2} \\
		& + \alpha_k^2 \mathbb{E}[\psi_k \indicator_{H_k}] \label{eqn:Rk_term3}.
\end{align}
\end{subequations}
The term in the right-hand side of (\ref{eqn:Rk_term1}) is 0 since $\mathbb{E}[Z_k] = 0$. For the term (\ref{eqn:Rk_term2}) we have
\begin{equation}
\begin{aligned}
	\mathbb{E}[\phi_k^2 \indicator_{H_k}] & = \mathbb{E}\left[ \indicator_{H_k} \left( Z_k^\mathrm{T} \nabla c(\mu_k) \right)^2 \right] \\ & \leq \mathbb{E}\left[ \indicator_{H_k} \norm{Z_k}^2 \norm{\nabla c(\mu_k)}^2 \right] \\
		& \leq \mathbb{E}\left[ \indicator_{E_{k+1}} \norm{Z_k}^2 \norm{\nabla c(\mu_k)}^2 \right] \\
		& \leq \sigma^2 K^2,
\end{aligned}
\end{equation}
where the first inequality is due to Cauchy-Schwarz, the second inequality is due to Lemma \ref{lem:Ek_Hk} and the third inequality is due to Lipschitz continuity of $c(\mu_k)$. Finally, for the term (\ref{eqn:Rk_term3}) we have
\begin{equation}
\begin{aligned}
	& \mathbb{E}\left[ \psi_k \indicator_{H_k} \right] \\ 
 & = \mathbb{E}\left[ \indicator_{H_k} 2L \left( \norm{\nabla c(\mu_k)}^2 + \norm{Z_k}^2 \right) \right] \\
 	& \leq 2L \left( K^2 + \sigma^2 \right).
\end{aligned}
\end{equation}
The term $R_{*} \alpha_k^2$ in (\ref{eqn:Rk_lemma}) follows immediately from the bounds for (\ref{eqn:Rk_term2}) and (\ref{eqn:Rk_term3}) derived above. The last step is to show how the term $\varepsilon \mathbb{P}\left[ \widetilde{H}_{k+1} \right]$ in (\ref{eqn:Rk_lemma}) is obtained from $R_k \indicator_{\widetilde{H}_{k+1}}$. Note that by definition we have $R_{k} > \varepsilon$ if $\widetilde{H}_{k+1}$ occurs. Hence
\begin{equation}
\mathbb{E}\left[ R_k \indicator_{\widetilde{H}_{k+1}} \right] \geq \varepsilon \mathbb{E}\left[ \indicator_{\widetilde{H}_{k+1}} \right] = \varepsilon \mathbb{P}\left[ \widetilde{H}_{k+1} \right]
\end{equation}
and the result follows.
$\blacksquare$
\end{pf}
With the above derivations completed, we now have the necessary ingredients in place to prove Theorem \ref{thm:stability}. First, recall that since $H_{k} \subseteq E_{k+1}$, we have $\mathbb{P}[E_{k+1}] \geq \mathbb{P}[H_{k}]$ and it is sufficient to prove that $\mathbb{P}\left[ H_{k} \right] \geq 1 - \delta$ for all $k = 1, 2, ...$. Denoting by $\widetilde{H}_i^c$ the complement of $\widetilde{H}_i$, we have
\begin{equation}
\begin{aligned}
	\mathbb{P}[ H_{k} ] = \mathbb{P}\left[ \bigcap_{i=0}^{k} \widetilde{H}_i^c \right]
		& = 1 - \mathbb{P}\left[ \bigcup_{i=0}^{k} \widetilde{H}_i \right] \\
		& = 1 - \sum_{i=0}^{k} \mathbb{P} \left[ \widetilde{H}_i \right],
\end{aligned}
\label{eqn:P_Hk}
\end{equation}
where the first equality can be intuitively understood as the statement that ``$H_k$ occurs if there is no large noise in each iteration, i.e., $\widetilde{H}_i$ does \emph{not} occur". The third equality in (\ref{eqn:P_Hk}) is due to disjointness of the events $\widetilde{H}_k$ for all $k = 0, 1, ...$. Our claim thus holds if $\sum_{i=0}^{\infty} \mathbb{P}\left[ \widetilde{H}_i \right] \leq \delta$, which we prove by first telescoping (\ref{eqn:Rk_lemma}):
\begin{equation}
	\mathbb{E}\left[ \widetilde{R}_k \right] \leq R_* \sum_{i=0}^{k-1} \alpha_i^2 - \varepsilon \sum_{i=1}^{k} \mathbb{P}\left[ \widetilde{H}_i \right].
\label{eqn:Rk_telescoped}
\end{equation}
The event $\widetilde{H}_{k+1}$ can be bounded by
\begin{equation}
\begin{aligned}
\mathbb{P}\left[ \widetilde{H}_{k+1} \right] & = \mathbb{P}[ H_{k} \setminus H_{k+1} ] = \mathbb{P}[ H_{k} \cap \{ R_k > \varepsilon \}] \\
	& = \mathbb{E}[ \indicator_{H_{k}} \times \indicator_{R_k > \varepsilon} ] \\
	& \leq \mathbb{E} [\indicator_{H_{k}} \times (R_k / \varepsilon) ] \\
	& = \mathbb{E}\left[ \widetilde{R}_k \right] / \varepsilon,
\end{aligned}
\label{eqn:P_Hk_bound}
\end{equation}
where the inequality follows from the fact that $R_k \geq 0$. Application of the bound given by (\ref{eqn:P_Hk_bound}) to the left-hand side of (\ref{eqn:Rk_telescoped}) gives
\begin{equation}
	\varepsilon \mathbb{P}\left[ \widetilde{H}_{k+1} \right] \leq R_* \sum_{i=0}^{k-1} \alpha_i^2 - \varepsilon \sum_{i=1}^{k} \mathbb{P}\left[ \widetilde{H}_i \right]
\end{equation}
and hence
\begin{equation}
	\sum_{i=1}^{k+1} \mathbb{P}\left[ \widetilde{H}_i \right]\leq \frac{R_*}{\varepsilon} \sum_{i=0}^{k-1} \alpha_i^2. 
	\label{eqn:sum_P_Hk_bound}
\end{equation}
Finally, we can verify that (\ref{eqn:Einfty_bound}) holds by combining (\ref{eqn:P_Hk}), (\ref{eqn:sum_P_Hk_bound}) and the step size bound (\ref{eqn:alpha_sum_bound}):
\begin{equation}
\begin{aligned}
	\mathbb{P}[E_\infty] \geq \mathbb{P}[H_\infty] & = 1 - \sum_{i=0}^{\infty} \mathbb{P}\left[ \widetilde{H}_i \right] \\
	& \geq 1 - \frac{R_*}{\varepsilon} \sum_{i=0}^{\infty} \alpha_i^2 \geq 1 - \delta.
\end{aligned}
\end{equation}

\section{Regularity of the cost function}

\begin{figure}[!htb]
\centering
\begin{tikzpicture}[auto, node distance=2cm,>=latex']
    \node [input, name=u1_input] 	(u1_input) {};
    \node [sum, below of=u1_input] (sum1) {$+$};
	\node [block, left of=sum1] 	(P) {$P$};
	\node [tmp, below of=sum1] 		(K_y1) {};
	\node [block, left of=K_y1] 	(K) {$K(\mu)$};
	\node [output, name=y1_output, below of=K_y1, node distance=1cm] (y1_output) {};
	\node [sum, left of=K, node distance=2cm] 			(sum2) {$+$};
	\node [input, name=u2_input, below of=sum2, node distance=1cm] (u2_input) {};
	\node [tmp, left of=P] 			(P_y2) {};
	\node [output, above of=P_y2, node distance=1cm]  (y2_output) {};
    
    \draw [->] (P) -- (sum1);
    \draw [->] (u1_input) -- node{$u_1$} (sum1);
    \draw [-] (sum1) -- (K_y1);
    \draw [->] (K_y1) -- (K);
    \draw [->] (K_y1) -- node{$y_1$} (y1_output);
    \draw [->] (K) -- (sum2);
    \draw [->] (u2_input) -- node{$u_2$} (sum2);
    \draw [-] (sum2) -- (P_y2);
    \draw [->] (P_y2) -- (P);
    \draw [->] (P_y2) -- node{$y_2$} (y2_output);
    \end{tikzpicture}
\caption{Decomposition of $G(\mu): (u_1, u_2) \rightarrow (y_1, y_2)$ into a large-scale $P$ and a small-scale parameter-dependent part $K(\mu)$.} \label{fig:PK_form}
\end{figure}
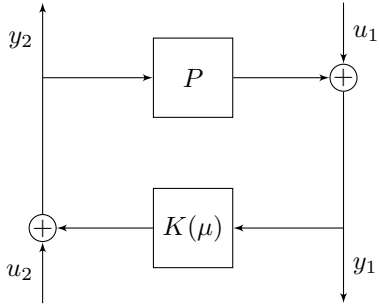

The regularity assumptions of the cost function (\ref{eqn:cost}) are given in Assumptions \ref{assum:c_continuity} and \ref{assum:c_sublevel_lipschitz}. These assumptions are required for the stability preservation result stated in Theorem \ref{thm:stability}. In this section, we provide some additional insight to understand under which conditions Assumption \ref{assum:c_continuity} is guaranteed to hold. Specifically, we consider the parameter-dependent system-to-be-optimized $G(\mu)$ as a feedback interconnection of a (large-scale) constant sub-system $P$ and a (small-scale) parameter-dependent sub-system $K(\mu)$ (Figure \ref{fig:PK_form}). 

\subsection{Auxiliary results}
The following lemma is from \cite[Theorem 1]{Zedek1965} and will be used in Section \ref{sct:continuity}.
\begin{lem} \label{lem:polynomial_roots_continuous}
Given a polynomial $p_n(s) = \sum_{k=0}^{n} a_k s^k$ with $a_n \neq 0$, an integer $m \geq n$ and a number $\varepsilon > 0$, there exist a number $\delta > 0$ such that whenever the $m+1$ complex numbers $b_k$, $0 \leq k \leq m$ satisfy the inequalities
\begin{equation}
\begin{aligned}
|b_k - a_k | & < \delta && \text{ for } 0 \leq k \leq n, \\
|b_k| & < \delta && \text{ for } n+1 \leq k \leq m,
\end{aligned}
\end{equation}
then the zeros $\beta_k$, $1 \leq k \leq m$, of the polynomial $q_m(s) = \sum_{k=0}^{m} b_k s^k$ can be labeled in such a way as to satisfy with respect to the zeros $\alpha_k$, $1 \leq k \leq n$, of $p_n(s)$ the inequalities
\begin{equation}
\begin{aligned}
|\beta_k - \alpha_k| & < \varepsilon  && \text{ for } 1 \leq k \leq n, \\
|\beta_k| & > \varepsilon^{-1} 	   && \text{ for } n+1 \leq k \leq m.
\end{aligned}
\end{equation}
\end{lem}

\subsection{Main result} \label{sct:continuity}
\begin{thm} \label{thm:cost_continuity}
Let $G(\mu): (u_1, u_2) \rightarrow (y_1, y_2)$ admit a feedback decomposition as depicted in Figure \ref{fig:PK_form}. Then, $c(\mu) = \frac{1}{2} \hnorm{G(\mu)}^2$ is continuous for all $\mu \in \muSpace$ if the following conditions are met:
\begin{enumerate}
	\item[(i)] The McMillan degree of $K(\mu)$ is fixed.
	\item[(ii)] $P(s)$ is a rational transfer function.
	\item[(iii)] Each element of the transfer matrix $K(\mu; s)$ is a rational transfer function for which the coefficients in the numerator and denominator are continuous functions of $\mu$.
	\item[(iv)] $G(\mu)$ is strictly proper.
\end{enumerate}
\end{thm}
\begin{pf}
Condition (i) implies that the McMillan degree (and hence the number of poles) of $G(\mu)$ does not depend on $\mu$ and its poles are equal to the poles of $(I - PK(\mu))^{-1}$ \cite[Theorem 30.1]{dahleh2004lectures}. A pole of $(I - PK(\mu))^{-1}$ is a zero of the polynomial $\mathrm{det}(I - P(s)K(s; \mu))$. From conditions (ii) and (iii) we know that the entries of $P(s)$ and $K(s; \mu)$ are polynomials in $s$ (with coefficients that depend continuously on $\mu$ in the case of $K(s; \mu)$). By extracting the least common denominator, they can be written as
\begin{align}
	P(s) & = \frac{1}{d_P(s)} \bar{P}(s), \\
	K(s; \mu) &= \frac{1}{d_K(s; \mu)} \bar{K}(s; \mu),
\end{align}
where $d_P(s)$ is a polynomial, $\bar{P}(s)$ is a matrix of polynomials, $d_K(s; \mu)$ a polynomial with $\mu$-dependent coefficients and $\bar{K}(s; \mu)$ a matrix of polynomials with $\mu$-dependent coefficients. Using this representation, the zeros of $\mathrm{det}(I - P(s)K(s; \mu))$ are equivalent to the zeros of the polynomial
\begin{equation}
	\phi(s; \mu) := \mathrm{det}\left( d_P(s) d_K(s; \mu) I - \bar{P}(s) \bar{K}(s; \mu) \right).
	\label{eqn:phi}
\end{equation}
Since $\phi(s; \mu)$ is a polynomial whose coefficients depend continuously on $\mu$, we can apply Lemma \ref{lem:polynomial_roots_continuous} to conclude that its roots (and hence the poles of $G(\mu)$) depend continuously on $\mu$. To show how continuous dependence of the poles on $\mu$ leads to continuity in the $\mathcal{H}_2$ norm, we split $\muSpace$ (the domain of $\mu$) into 3 parts:
\begin{enumerate}
	\item $\mu \in \mathrm{Interior}(\muSpaceStab)$ (all poles are in the open left-half plane).
	\item $\mu \in \muSpace \setminus \muSpaceStab$ (at least one pole is in the closed right-half plane and cannot be stabilized by an arbitrarily small perturbation in $\mu$).
	\item $\mu \in \mathrm{Boundary}(\muSpaceStab)$ (all poles are either in the open left-half plane or they are on the stability boundary ($i\mathbb{R}$) and can be stabilized by an arbitrarily small perturbation in $\mu$).
\end{enumerate}

For case 1), by continuity of the poles of $G(\mu)$, we have that for each $\mu \in \mathrm{Interior}(\muSpaceStab)$, there exists a $\delta > 0$ such that $G(\mu +  x)$ is asymptotically stable for all $\norm{x} < \delta$. Stability and condition (iv) together imply boundedness of (\ref{eqn:h2inner}) and continuity of $c(\mu)$ then follows because it is equivalent to continuity of the integral (\ref{eqn:h2inner}).

For case 2), again by continuity of the poles of $G(\mu)$, we have that for each $\mu \in \mathrm{Interior}(\muSpace \setminus \muSpaceStab)$, there exists a $\delta > 0$ such that $G(\mu + x)$ is unstable for all $\norm{x} < \delta$. Continuity follows since $c(\mu) = c(\mu + x) = \infty$.

Case 3) is more involved. Since $\mu \in \mathrm{Boundary}(\muSpace \setminus \muSpaceStab)$, we have $c(\mu) = \infty$. For continuity, we need to show that for each (arbitrarily small) $\varepsilon > 0$ there exists a $\delta > 0$ such that $c(\mu + x) > \varepsilon^{-1}$ for all $\norm{x} < \delta$. This holds trivially for perturbations $x$ that do not stabilize the poles (since then $c(\mu + x) = \infty$), so we only need to consider perturbations $x$ that move all unstable poles from the imaginary line into the open left half plane. Fix such an arbitrary perturbation $x$ and introduce a variable $\alpha \in [0, 1]$ that sweeps the parameter according to $\mu + \alpha x$. Continuity of $G(\mu)$ is thus equivalent to showing that $\lim_{\alpha \rightarrow 0} G(\mu + \alpha x) = \infty$.
Let $\{ p_i(\alpha) \}_{i=1}^{n_p}$ denote the subset of the poles of $G(\mu)$ that are on the stability boundary for $\alpha = 0$ and, by our choice of a stabilizing perturbation $x$, in the open left half plane for $\alpha \in (0, 1]$.
Assume that these poles are distinct for $\alpha \in (0, 1]$, so that we can decompose $G(\mu + \alpha x)$ as follows:
\begin{equation}
\begin{aligned}
	G(s; \mu + \alpha x) = \underbrace{\sum_{i=1}^{n_p} \frac{c_i(\alpha) b_i(\alpha)^\mathrm{T}}{s - p_i(\alpha)}}_{P(s; \alpha)} + S(s; \alpha).
\end{aligned}
\end{equation}
(If, on the other hand, not all poles are distinct, then we can take powers of the denominators in $P(s; \alpha)$ according to the multiplicity of these poles and continue with the proof in the same manner).
By condition (iv) in the statement of the theorem, $P(\alpha)$ and $S(\alpha)$ are both strictly proper.
Based on this decomposition, we have the following lower bound:
\begin{equation}
\begin{aligned}
\hnorm{G(\mu + \alpha x)}^2 \geq \hnorm{P(\alpha)}^2 + 2 \hinner{ P(\alpha) }{ S(\alpha) }.
\end{aligned}
\label{eqn:PS_decomposition}
\end{equation}
By showing that $\hnorm{P(\alpha)}^2$ is unbounded and $\hinner{P(\alpha)}{S(\alpha)}$ is finite as $\alpha \rightarrow 0$, we will prove the theorem. For the first term of the right-hand side of bound (\ref{eqn:PS_decomposition}) we write using the residue-decomposition of the $\mathcal{H}_2$ norm (see, e.g., \cite[Lemma 2.1.4]{Antoulas2020}):
\begin{equation}
\begin{aligned}
	\hnorm{P(\alpha)}^2 & = \sum_{i=1}^{n_p} c_i(\alpha)^\mathrm{T} P(-p_i(\alpha); \alpha) b_i(\alpha) \\
		& = \sum_{i=1}^{n_p} \sum_{j=1}^{n_p} \frac{ c_i(\alpha)^\mathrm{T} c_j(\alpha) b_j(\alpha)^\mathrm{T} b_i(\alpha) }{-\overline{p_i}(\alpha) - p_j(\alpha)}.
\end{aligned}
\label{eqn:P_decomp}
\end{equation}
The residues $c_i(\alpha), b_i(\alpha)$ are non-zero for $\alpha \in [0, 1]$ since the McMillan degree of $G(\mu)$ is fixed (condition (i) in the theorem). Hence, for $i = j$ in the last term of (\ref{eqn:P_decomp}) the numerator is positive while the denominator goes to 0 (since $-\overline{p_i}(\alpha) - p_i(\alpha) = -2 \mathrm{Re}(p_i(\alpha) \rightarrow 0$ as $\alpha \rightarrow 0$). For $i \neq j$ the denominator is non-zero since the poles $p_i(\alpha)$, $i = 1, ..., n_p$, are distinct for $\alpha \in (0, 1]$.
Therefore, we have for the limit that
\begin{equation}
	\lim_{\alpha \rightarrow 0} \hnorm{P(\alpha)}^2 = \lim_{\alpha \rightarrow 0} \sum_{i=1}^{n_p} \frac{ \norm{c_i(\alpha)} \norm{b_i(\alpha)}}{-2 \mathrm{Re}(p_i(\alpha))} = \infty,
\end{equation}
To show boundedness of the second term of bound (\ref{eqn:PS_decomposition}) we write
\begin{equation}
	\hinner{P(\alpha)}{S(\alpha)} = \sum_{i=1}^{n_p} c_i(\alpha)^\mathrm{T} S(-\overline{p_i}(\alpha)) b_i(\alpha).
	\label{eqn:PS}
\end{equation}
By definition, the poles of $S(\alpha)$ do not approach stability boundary ($i\mathbb{R}$) as $\alpha \rightarrow 0$ (namely, such poles are in $P(\alpha)$). Thus, all terms on the right-hand side of (\ref{eqn:PS}) are finite and hence
\begin{equation}
	\lim_{\alpha \rightarrow 0} | \hinner{P(\alpha}{S(\alpha)} | = \left| \sum_{i=1}^{n_p} c_i(0)^\mathrm{T} S(-\overline{p_i}(0)) b_i(0) \right| < \infty.
\end{equation}
Concluding, continuity of $G(\mu)$ is shown for all $\mu \in \muSpace$.
$\blacksquare$
\end{pf}

\end{appendix}

\end{document}